\documentclass[11pt,twoside]{article}
\usepackage[utf8]{inputenc}
\usepackage[english]{babel}
\usepackage{fancyhdr}
\usepackage{amssymb,amsmath,amsthm}
\usepackage{blindtext}
\usepackage{hyperref}
\usepackage{caption}
\usepackage[square,numbers]{natbib}

\usepackage{mathtools}
\usepackage{enumerate}
\usepackage{amssymb,amsmath}
\usepackage{algorithm}
\usepackage{algpseudocode}
\usepackage{multirow}
\usepackage{placeins}
\usepackage{bigstrut}

\newcommand{\R}{\mathbb R} 
 \DeclarePairedDelimiter{\norm}{\lVert}{\rVert} 
 \DeclarePairedDelimiter{\abs}{\lvert}{\rvert} 
 \DeclareMathOperator*{\argmin}{Argmin} 
 \DeclareMathOperator*{\argmax}{Argmax} 

\newtheorem{assumption}{Assumption}
\newtheorem{definition}{Definition}
\newtheorem{lemma}{Lemma}
\newtheorem{proposition}{Proposition}
\newtheorem{theorem}{Theorem}

\newtheorem{remark}{Remark}

\newcommand{\ASSIMPL}{\texttt{AS-SIM\-PLEX}}
\newcommand{\FW}{\texttt{FW}}
\newcommand{\AFW}{\texttt{AFW}}
\newcommand{\PG}{\texttt{PG}}
\newcommand{\ASFW}{\texttt{AS-FW}}
\newcommand{\ASAFW}{\texttt{AS-AFW}}
\newcommand{\ASPG}{\texttt{AS-PG}}
\newcommand{\eps}{\epsilon}

\hypersetup{
    colorlinks=false,
    urlcolor=black,
    pdfborder={0 0 0}
}

\begin{document}
\thispagestyle{plain}

\setcounter{page}{1}

{\centering
{\LARGE \bfseries An Active-Set Algorithmic Framework for Non-Convex Optimization Problems over the Simplex}

\bigskip\bigskip
Andrea Cristofari$^*$, Marianna De Santis$^\dag$, Stefano Lucidi$^\dag$, Francesco Rinaldi$^*$
\bigskip

}

\begin{center}
\small{\noindent$^*$Dipartimento di Matematica ``Tullio Levi-Civita'' \\
Universit\`a di Padova \\
Via Trieste, 63, 35121 Padua, Italy \\
E-mail: \texttt{andrea.cristofari@unipd.it}, \texttt{rinaldi@math.unipd.it} \\
\bigskip
$^\dag$Dipartimento di Ingegneria Informatica, Automatica e Gestionale \\
Sapienza Universit\`{a} di Roma \\
Via Ariosto, 25, 00185 Rome, Italy \\
E-mail: \texttt{mdesantis@diag.uniroma1.it}, \texttt{lucidi@diag.uniroma1.it} \\
}
\end{center}

\bigskip\par\bigskip\par
\noindent \textbf{Abstract.}
In this paper, we describe a new active-set algorithmic framework for minimizing
a non-convex function over the unit simplex.  At each iteration, the method makes use of a
rule for identifying active variables (i.e., variables that are zero at a stationary point) and specific
directions (that we name \textit{active-set gradient related} directions) satisfying  a new ``nonorthogonality'' type of condition.
We prove global convergence to stationary points when using an Armijo line search in the given framework.
We further describe three different examples of active-set gradient related directions that
guarantee linear convergence rate (under suitable assumptions).
Finally, we report numerical experiments showing the effectiveness of the approach.

\bigskip\par
\noindent \textbf{Keywords.} Active-set methods. Unit simplex. Non-convex optimization. Large-scale optimization.

\bigskip\par
\noindent \textbf{MSC2000 subject classifications.} 65K05. 90C06. 90C30.

\section{Introduction}\label{sec:intr_simplex}
Many real-world applications can be modeled as optimization problems over structured feasible sets.
In particular, the problem of minimizing a function over a simple polytope
(such as the unit simplex)
arises in different fields like, e.g., machine learning, statistics and economics.
Examples of relevant applications include  training of support vector machines, boosting (Adaboost),
convex approximation in $\ell_p$, mixture density estimation, finding maximum stable
sets (maximum cliques) in graphs, portfolio optimization and population dynamics problems
(see, e.g.,~\cite{bomze1997evolution,bomze:1999,clarkson2010coresets,deKlerk:2008} and references therein).

The problem we address here can be stated as follows:
\begin{equation}\label{simplex_prob}
      \begin{array}{l}
      \displaystyle\min_{x\in \Delta}~  f(x) \\[1ex]
      \end{array}
\end{equation}
where $\Delta=\{x\in \R^n:\ e^T x = 1,~x \ge 0\}$ is the unit simplex, $e\in \R^n$
is the vector of all ones, $f \colon \R^n \rightarrow \R$ is continuously differentiable and its gradient
$\nabla f(x)$ is Lipschitz continuous over the feasible set, with constant $L>0$.

Note that minimizing an objective function $h(x)$ over a polytope $P$ can be recast as problem~\eqref{simplex_prob}.
Indeed, since any point $x\in P$ can be expressed as a convex combination of the columns
of \mbox{$V = \begin{bmatrix} v_1 & \ldots & v_m \end{bmatrix} \in \R^{n \times m},$}
with $v_1,\dots,~v_m$ vertices of $P$, the problem $\min\{h(x): \, x\in P\}$ can be rewritten as
\hbox{$\min\{h(Vy) : \, e^T y = 1, \, y \ge 0\}$}. Thus, each variable $y_i$ represents
the weight of the $i$th vertex in the convex combination.

In many different contexts, problems can have very sparse solutions (i.e., solutions with many zero components).
Hence, developing methods that allow to quickly identify the set of zero components in the optimal solution
is getting crucial to guarantee relevant savings in terms of CPU time.
In our problem, estimating the set of zero components in the optimal solution, or in a stationary point when the objective function is non-convex,
coincides with estimating the set of \textit{active} (or binding) inequality constraints.
This set of active constraints is often referred to as \textit{active set} and the so called \textit{active-set methods}
are characterized by computing, at each iteration, an estimate of the binding constraints which is iteratively updated.
Usually, only a single active constraint is added to or deleted
from the estimated active set at each iteration (see, e.g.,~\cite{nocedal:2006} and references therein). However, when dealing with simple constraints,
more sophisticated active-set methods can be used, which can add to or delete
from the current estimated active set more than one constraint at each
iteration, and eventually find the active set in a finite number of steps if certain conditions hold.
In particular, several active-set algorithms
are based on the idea of combining an ``identification'' step (to identify the active constraints)
with a minimization step in a reduced space
(see, e.g.,~\cite{bertsekas:1982,birgin2002large,bras:2017,cristofari:2017,diserafino:2017,facchinei1998accurate,facchinei1998active,hager:2006,hager:2016a,hager:2016b,more1989algorithms} and references therein).

Here, we propose an active-set algorithmic framework for solving problem~\eqref{simplex_prob}, where $f(x)$ is a possibly non-convex objective function.
Hereinafter, we call \textit{active variables} those variables that are equal to zero at a stationary point.
Indeed, those variables correspond to the active constraints when the feasible set is the unit simplex.

In the fist part of the paper, we describe an active-set estimate to identify the active variables by extending some specific strategies proposed in the
contexts of box-constrained problems (see~\cite{buchheim:2016,cristofari:2017,desantis:2012,desantis:2016})
to the case of unit simplex. The main features of our active-set strategy are essentially two:
\begin{enumerate}
 \item it does not only focus on the zero variables and keep them fixed, but rather tries to
 quickly identify as many active variables as possible (including nonzero variables) at a given point;
 \item it gives a significant reduction in the objective function, while guaranteeing feasibility,
  when setting to zero those variables estimated active (and moving a suitably chosen variable estimated nonactive).
\end{enumerate}
These features enable us to easily use this active-set strategy into any globally convergent algorithm, potentially improving its efficiency without
compromising its theoretical properties.



In the second part of the paper, inspired by the ``nonorthogonality'' type of condition described in~\cite{bertsekas1999nonlinear},
we define a new class of directions, named \textit{active-set gradient related},
that we combine with the active-set strategy described above to devise a two-step algorithmic framework.
In the first step, our method sets to zero the estimated active variables and suitably moves an estimated nonactive variable,
producing a new feasible point with a smaller objective value.
In the second step, an active-set gradient related direction, combined with a suitable line search,
is used in the subspace of the estimated nonactive variables to generate the next iterate.
We prove convergence of our framework to stationary points of problem~\eqref{simplex_prob} using Armijo line search.
We would like to highlight that, since at each iteration we move to zero some of the variables and then approximately
optimize in a subspace, guaranteeing convergence is not a straightforward task and requires a thorough theoretical analysis.

We further give three specific examples of active-set gradient related directions that can be
used to implement our algorithm in practice. More specifically, we consider Frank-Wolfe~\cite{frank_wolfe:1956}, away-step Frank-Wolfe~\cite{wolfe:1970} and
projected gradient directions (see, e.g., \cite{bertsekas1999nonlinear} and references therein).
We then prove linear convergence rate of the framework (under suitable assumptions) when using these directions.

In the final part of the paper, we report numerical results on both convex and non-convex instances. The results seem to indicate that
our active-set algorithm is very efficient when dealing with sparse optimization problems.

The paper is organized as follows. In Section~\ref{sec:notation}, we report the notation and useful preliminary results.
In Section~\ref{sec:act_set_simplex}, we describe in depth our active-set estimate and the related theoretical properties.
In Section~\ref{sec:algorithm_simplex}, we present our algorithmic framework and
carry out the convergence analysis. We also give three specific examples of active-set gradient related directions that can be used in the framework.
In Section~\ref{sec:convrate}, we analyze the convergence rate of the method when using those active-set gradient related directions.
In Section~\ref{sec:res_simplex}, we report our numerical experience.
Finally, in Section~\ref{sec:conclusions_simplex}, we draw some conclusions.

\section{Notation and Preliminaries}\label{sec:notation}
Throughout the paper, we indicate with $\|\cdot\|$ the Euclidean norm.
Given a vector $v\in \R^n$ and an index set $I \subseteq \{1,\ldots,n\}$, we denote with $v_I$ the subvector with components $v_i$, $i\in I$.
We indicate with $e_i$ the $i$th unit vector.
Given $x \in \R^n$ and a non-empty closed convex set $S \subseteq \R^n$, we denote by $P(x)_S$ the projection of $x$ on $S$.
Finally, given a subset $I \subseteq \{1,\ldots,n\}$, we denote $\Delta_I$
the subset of points of $\Delta$ with $x_i = 0$ for all $i \notin I$, i.e.,
\[
\Delta_I := \{x \in \Delta \colon x_i =0, \; \forall i \notin I\}.
\]

\begin{definition}
A stationary point $x^*$ of problem~\eqref{simplex_prob} is a feasible point that satisfies
the following first-order necessary optimality conditions:
\begin{gather}
\nabla f(x^*) - \lambda^* e - \mu^* = 0, \label{opt1}\\
(\mu^*)^T x^* = 0, \label{sc2}\\
\mu^* \ge 0. \label{non2}
\end{gather}
where $\lambda^* \in \R$ and $\mu^* \in \R^n$ are the KKT multipliers.
\end{definition}

It is easy to verify that conditions~\eqref{opt1}--\eqref{non2} are equivalent to the following:
\[
\nabla_i f(x^*)
\begin{cases}
\ge \lambda^*, \quad & x^*_i = 0, \\
= \lambda^*, \quad & x^*_i > 0,
\end{cases}
\qquad \qquad i = 1,\ldots,n.
\]
\begin{remark}\label{remark:kkt}
Note that a feasible point $x^*$ of problem~\eqref{simplex_prob} is stationary if and only if $\nabla f(x^*)^T(e_i - x^*) \ge 0$ for all $i = 1,\ldots,n$, or equivalently,
\[
\max\{0, -\nabla f(x^*)^T(e_i - x^*)\} = 0, \quad i = 1,\ldots,n.
\]
\end{remark}

\section{Active-Set Estimate}\label{sec:act_set_simplex}
Given a stationary point $x^*$ of problem~\eqref{simplex_prob}, the active set can be defined as the set of inequality constraints binding at $x^*$.
As mentioned above, since there is a one-to-one correspondence between inequality constraints and variables, we equivalently define as active set
the set of zero components at $x^*$.
\begin{definition}\label{def:activeset}
Let $x^*\in \R^n$ be a stationary point of problem~\eqref{simplex_prob}.
We define as active set the following set:
\begin{equation}\label{AS}
\bar A (x^*) = \bigl \{i\in \{1,\ldots,n\} \colon x^*_i=0 \bigr\}.
\end{equation}
We further define the nonactive set $\bar N (x^*)$ as the complement of $\bar A(x^*)$:
\begin{equation}\label{NAS}
\bar N (x^*) = \{1,\ldots,n\} \setminus \bar A (x^*) = \bigl\{i\in \{1,\ldots, n\} \colon x^*_i > 0 \bigr\}.
\end{equation}
\end{definition}
Now, we describe how, at any feasible point $x$, we estimate the active set.
Following the approach proposed in~\cite{dipillo:1984,facchinei:1995},
we use a strategy that requires proper approximation of the KKT multipliers by means of
the so called \emph{multiplier functions}.
To compute these multiplier functions, let $(\lambda^*,\mu^*)$ be the KKT multipliers associated to a given stationary point $x^*$.
By~\eqref{opt1}, we have
\begin{equation*}
\mu^*= \nabla f(x^*)-\lambda^* e,
\end{equation*}
then, multiplying by $x^*$ and taking into account complementarity condition~\eqref{sc2}, we get
\begin{equation*}
0=(\mu^*)^Tx^*= (\nabla f(x^*)-\lambda^* e)^Tx^*.
\end{equation*}
From the feasibility of $x^*$, we obtain the following expressions for the multipliers:
\begin{gather}
\lambda^*= \nabla f(x^*)^Tx^*,\quad \mu^*= \nabla f(x^*)- \lambda^* e, \nonumber
\end{gather}
so that we can introduce the following  multiplier functions:
\begin{gather}
\lambda(x) = \nabla f(x)^T x, \label{lambdamumult} \\
\mu_i(x) = \nabla_i f(x)- \lambda(x), \quad i = 1,\ldots,n. \label{mumumult}
\end{gather}
Now, we can define our estimate of the active set.
\begin{definition}\label{def:activesetest}
Let $x\in \R^n$ be a feasible point of problem~\eqref{simplex_prob}.
We define the active-set estimate $A(x)$ and the nonactive-set estimate $N(x)$ as
\begin{gather}
A(x) = \{i \colon x_i \le \epsilon \mu_i(x) \bigr\}  = \{i \colon x_i \le \epsilon \nabla f(x)^T(e_i - x) \}, \label{as1} \\
N(x) = \{i \colon x_i > \epsilon \mu_i(x)\} = \{i \colon x_i > \epsilon \nabla f(x)^T(e_i - x)\},  \label{nas1}
\end{gather}
where $\epsilon$ is a positive scalar.
\end{definition}

By adapting Theorem~2.1 in~\cite{facchinei:1995},
we can ensures that, in a neighborhood of a stationary point $x^*$, all the estimated active variables are active at $x^*$ and
include all active variables at $x^*$ satisfying strict complementarity.
We state this result in the following theorem.
\begin{theorem}\label{th:estim}
If $(x^*,\lambda^*,\mu^*)$ satisfies KKT conditions for problem~\eqref{simplex_prob}, then there exists
an open ball $\mathcal B(x^*,\rho)$ with center $x^*$ and radius $\rho > 0$
such that, for each $x \in \mathcal B (x^*,\rho)$ we have
\begin{equation*}
\{i \colon x_i^*=0, \, \mu_i(x^*)>0\} \subseteq A(x) \subseteq \bar A (x^*).
\end{equation*}
Furthermore, if strict complementarity holds, then
\begin{equation*}
\{i \colon x_i^*=0, \, \mu_i(x^*)>0\} = A(x) = \bar A (x^*),
\end{equation*}
for each $x \in \mathcal B (x^*,\rho)$.
\end{theorem}

\subsection{A Global Property of the Active-Set Estimate}\label{sub:act_set}
Here, we analyze a global property of our active-set estimate.
In particular, we show how, given a point $x\in \R^n$ feasible for problem~\eqref{simplex_prob}, we can obtain a sufficient decrease in the objective function
by setting the estimated active variables to zero.
In the context of box constrained problems, this property has been already proved in~\cite{buchheim:2016,cristofari:2017,desantis:2012,desantis:2016}.
When dealing with problems over the standard simplex, we have to handle the presence of an equality constraint:
in order to maintain feasibility, we need to update at least one estimated
nonactive variable, so that all variables sum up to $1$. The next proposition gives us a hint on how to choose
the estimated nonactive variable that will be updated when setting to zero the active variables.
\begin{proposition}\label{prop0}
Let $x \in \R^n$ be a feasible non-stationary point of problem~\eqref{simplex_prob} and define
\begin{equation}\label{ind_set_j}
J(x) = \Bigl\{ j \colon j \in \displaystyle\argmin_{i=1,\ldots,n}\bigl\{\nabla_i f(x)\bigl\} \Bigr\}.
\end{equation}

Then,
$J(x) \subseteq N(x)$.
\end{proposition}
\begin{proof}
Since $x$ is non-stationary, we have $\abs{J(x)} < n$. Moreover, an index $i$ must exist such that $x_i > 0$ and $\nabla_i f(x) > \nabla_j f(x)$, $j \in J(x)$.
It follows that
\begin{equation*}
\nabla f(x)^Tx>\nabla_j f(x)e^Tx=\nabla_j f(x).
\end{equation*}
Now, we can choose any index $j \in J(x)$ and set $\nu = j$. Recalling the definition of multipliers~\eqref{lambdamumult}--\eqref{mumumult}, we obtain
\begin{equation*}
\mu_{\nu}(x) = \nabla_{\nu} f(x)-\lambda(x) = \nabla_{\nu} f(x)-\nabla f(x)^Tx < \nabla_{\nu} f(x)-\nabla_{\nu} f(x) = 0\le x_{\nu}.
\end{equation*}
Since $x_{\nu}\ge0$ and $\mu_\nu(x)<0$, we have that $x_{\nu} > \eps \mu_{\nu}(x)$, and then $\nu \in N(x)$.
\end{proof}

\begin{remark}\label{rem:N}
Proposition~\ref{prop0} implies that for every feasible non-stationary point $x$, the estimated nonactive set $N(x)$ is non-empty.
\end{remark}

The main result of this section, reported in Proposition~\ref{prop_as}, shows that it is possible to get a significant decrease in the objective function
by setting to zero the estimated active variables and suitably updating a variable chosen in the set defined in Proposition~\ref{prop0}.
To obtain this result, we first need an assumption on the parameter $\eps$ appearing in Definition~\ref{def:activesetest}.
\begin{assumption}\label{ass:eps}
Assume that the parameter $\eps$ appearing in the estimates~\eqref{as1}--\eqref{nas1} satisfies the following conditions:
\[
0 < \epsilon \le \frac 2{nL(2C+1)},
\]
where $C > 0$ is a given constant.
\end{assumption}

\begin{proposition}\label{prop_as}
Let Assumption~\ref{ass:eps} hold.
Given a feasible non-stationary point $x$  of problem~\eqref{simplex_prob}, let $j \in J(x)$
and $I = \{1,\ldots,n\} \setminus \{j\}$. Let ${\hat A}(x)$ be a set of indices such that $\hat A(x) \subseteq A(x).$
Let $\tilde x$ be the feasible point defined as follows:
\begin{equation*}
\tilde x_{\hat A(x)} = 0; \quad \tilde x_{I\setminus \hat A(x)} = x_{I\setminus \hat A(x)}; \quad \tilde x_j = x_j + \displaystyle{\sum_{i \in \hat A(x)} x_i}.
\end{equation*}
Then,
\begin{equation*}
f(\tilde x)-f(x)\le - C L \norm{\tilde x-x}^2,
\end{equation*}
where $C > 0$ is the constant appearing in Assumption~\ref{ass:eps}.
\end{proposition}
\begin{proof}
Define
\begin{equation}\label{eq:hatAp}
\hat A^+ = \hat A(x) \cap \{i \colon x_i > 0\}.
\end{equation}
Using known results on functions with Lipschitz continuous gradient~\cite{nesterov:2013}, we can write
\[
f(\tilde x) \le f(x) + \nabla f(x)^T (\tilde x-x) +  \frac L 2 \norm{\tilde x-x}^2.
\]
Adding and subtracting $CL\norm{\tilde x-x}^2$ in the right-hand side of the above inequality, we get
\begin{equation}\label{lips}
f(\tilde x) \le f(x) + \nabla f(x)^T (\tilde x-x) + \frac {L(2C+1)} 2 \norm{\tilde x-x}^2 - CL \norm{\tilde x-x}^2.
\end{equation}
In order to prove the proposition, we need to show that
\begin{equation}\label{eqThesis}
\nabla f(x)^T (\tilde x-x) + \frac {L(2C+1)} 2 \|\tilde x-x\|^2 \le 0.
\end{equation}
From the definition of $\tilde x$, we get
\begin{equation}\label{normdiff}
\begin{split}
\norm{\tilde x-x}^2 & = \sum_{i\in \hat A^+} (x_i)^2 + \Bigg(\sum\limits_{i\in \hat A^+} x_i \Bigg)^2 \le \sum_{i\in \hat A^+} (x_i)^2
+\abs{\hat A^+}\sum_{i\in \hat A^+} (x_i)^2 \\
& = (\abs{\hat A^+}+1) x_{\hat A^+}^T x_{\hat A^+}
\end{split}
\end{equation}
and
\begin{equation}\label{ineqGrad}
\begin{split}
 \nabla f(x)^T (\tilde x-x) & = -\nabla_{\hat A^+}f(x)^T x_{\hat A^+} + \nabla_j f(x)\sum_{i\in \hat A^+} x_i \\
 & = x_{\hat A^+}^T \Bigl(\nabla_j f(x) e_{\hat A^+} - \nabla_{\hat A^+}f(x)\Bigr).
\end{split}
\end{equation}
From the definition of the index $j$, we have that $\nabla_i f(x)\ge \nabla_j f(x)$ for all $i \in \{1,\ldots,n\}$.
Therefore, we can write
\begin{equation}\label{ineqSimp}
 \sum_{i=1}^n \nabla_i f(x) x_i \ge \sum_{i=1}^n \nabla_j f(x) x_i = \nabla_j f(x) \sum_{i=1}^n  x_i =\nabla_j f(x).
\end{equation}
Recalling the active-set estimate and using~\eqref{ineqSimp}, we have that
\[
x_i \le \epsilon \Bigl(\nabla_i f(x)-\sum_{i=1}^n \nabla_i f(x) x_i\Bigr) \le \epsilon \bigl(\nabla_i f(x)- \nabla_j f(x)\bigr), \quad \forall i\in \hat A^+,
\]
so that, by~\eqref{normdiff}, we can write
\begin{equation}\label{normdiff2}
\norm{\tilde x-x}^2 \le \eps (\abs{\hat A^+} + 1) \, x_{\hat A^+}^T \Bigl( \nabla_{\hat A^+} f(x) - \nabla_j f(x) e_{\hat A^+} \Bigr).
\end{equation}
From~\eqref{ineqGrad} and~\eqref{normdiff2}, we get
 \begin{equation*}
 \begin{split}
 \nabla f(x)^T (\tilde x-x) + & \,\frac {L(2C+1)} 2 \norm{\tilde x-x}^2 \le
 x_{\hat A^+}^T \bigl[ \nabla_j f(x)e_{\hat A^+} - \nabla_{\hat A^+} f(x) \bigr] + \\[1ex]
 & \quad + \frac {L(2C+1)} 2 (\abs{\hat A^+}+1)\eps \, x_{\hat A^+}^T \Bigl( \nabla_{\hat A^+} f(x) - \nabla_j f(x) e_{\hat A^+} \Bigr) \\[1ex]
 & = \biggl( \frac {L(2C+1)} 2 (\abs{\hat A^+}+1)\eps - 1 \biggr) \, x_{\hat A^+}^T \Bigl( \nabla_{\hat A^+} f(x) - \nabla_j f(x) e_{\hat A^+} \Bigr) \\[1ex]
 & \le \biggl( \frac {L(2C+1)} 2 n \eps - 1\biggr) \, x_{\hat A^+}^T \Bigl( \nabla_{\hat A^+} f(x) - \nabla_j f(x) e_{\hat A^+} \Bigr),
 \end{split}
 \end{equation*}
where the last inequality follows from the non-negativity of \[x_{\hat A^+}^T \Bigl( \nabla_{\hat A^+} f(x) - \nabla_j f(x) e_{\hat A^+} \Bigr)\]
(implied by~\eqref{normdiff2}) and from the fact that $\abs{\hat A^+}+1 \le n$ (implied by Proposition~\ref{prop0}, see Remark~\ref{rem:N}).
Then,~\eqref{eqThesis} follows from the assumption we made on~$\epsilon$.
\end{proof}

\begin{remark}
In Assumption~\ref{ass:eps}, the upper bound of $\epsilon$ depends on $n$.
Actually, Proposition~\ref{prop_as} still holds by replacing the constant $n$ by $\abs{\hat A^+}+1$ in the upper bound of $\epsilon$,
where $\hat A^+$ is defined as in~\eqref{eq:hatAp}.
This follows from the fact that $n$ is only used to upper bound $\abs{\hat A^+}+1$ in the proof of Proposition~\ref{prop_as}.
Note that, in general, $\hat A^+$ might be considerably smaller than $n$, but it depends on both the specific point $x$ and $\epsilon$ itself.
So, for the sake of simplicity, in Assumption~\ref{ass:eps} we use the constant $n$, even if all the theoretical results of the paper would hold by using $\abs{\hat A^+}+1$ instead.
\end{remark}

\begin{remark}
From Assumption~\ref{ass:eps} and Proposition~\ref{prop_as}, we see that
there is a trade-off, depending on the constant $C$, between the magnitude of the upper bound of $\epsilon$ and the decrease in the objective function
guaranteed by Proposition~\ref{prop_as}.
Namely, for small values of $C$, large values of $\epsilon$ can be used, and then, from~\eqref{AS}, a major number of variables might be estimated active.
But the corresponding decrease in the objective function might be small. Vice versa, for large values of $C$, we have the opposite situation.
\end{remark}

The property described in the above proposition is crucial for the analysis of the algorithm framework that we carry out in the next section.
We remark that there is no way to get the same result from~\cite{dipillo:1984,facchinei:1995}, where a similar approach is used to estimate the active set.
Indeed, in those papers the authors deal with non-linear inequality constraints and there is no such a result like Proposition~\ref{prop_as}.
As a consequence, they cannot get the same algorithmic framework we describe in the present paper.

\section{An Active-Set Algorithmic Framework for Minimization over the Simplex}\label{sec:algorithm_simplex}
In this section, we describe in depth an algorithmic framework that embeds the active-set estimate explained in the previous section.
The framework performs two different steps at each iteration:
the first one to update the estimated active variables, and the second one to update the estimated nonactive variables.
The aim is to exploit as much as possible the properties of our estimate:
first, the ability to identify those active variables satisfying
strict complementarity after a sufficiently large number of iterations,
according to the results in Theorem~\ref{th:estim};
second, the ability to get a decrease in the objective function when moving the variables according to Proposition~\ref{prop_as}.

In particular, let $x^k$ be the point given at the beginning of a generic iteration $k$.
In the first step, we compute the active and nonactive-set estimates $A(x^k)$, $N(x^k)$,
and we generate the new feasible point $\tilde x^k$ as indicated in Proposition~\ref{prop_as}:
we set $\tilde x_{A(x^k)}$ to zero and we update $\tilde x^k_j$, with $j \in J(x^k)$ (all the other variables are not moved).
Then, in the second step, we compute a search direction $d_{N(x^k)}^k$ in the subspace of the estimated
nonactive variables, and we eventually perform a line search to get a new iterate $x^{k+1}$ (the computation of the search direction and the line search are described below).

The detailed scheme of our algorithmic framework, that we name \ASSIMPL, is reported in Algorithm~\ref{alg:AS_SIMPL}.

\begin{algorithm}[h!]
\caption{\texttt{Active-Set algorithmic framework for minimization over the simplex (\ASSIMPL)}}
\label{alg:AS_SIMPL}
\begin{algorithmic}
{\small
\par\vspace*{0.1cm}
\item[]$\,\,\,1$\hspace*{0.1truecm} Choose a feasible point $x^0$
\item[]$\,\,\,2$\hspace*{0.1truecm} For $k=0,1,\ldots$
\item[]$\,\,\,3$\hspace*{0.9truecm} If $x^k$ is a stationary point, then STOP
\item[]$\,\,\,4$\hspace*{0.9truecm} Compute $A^k = A(x^k)$ and $N^k = N(x^k)$
\item[]$\,\,\,5$\hspace*{0.9truecm} Compute $J^k = J(x^k)$, choose $j \in J^k$ and define $\tilde N^k = N^k\setminus\{j\}$
\item[]$\,\,\,6$\hspace*{0.9truecm} Set $\tilde x^k_{A^k} = 0\,$,
                                  $\,\tilde x^k_{\tilde N^k} = x^k_{\tilde N^k}\,$ and
                                  $\displaystyle{\,\tilde x^k_j = x^k_j+\sum_{h\in A^k} x^k_h}$
\item[]$\,\,\,7$\hspace*{0.9truecm} Compute a feasible direction $d^k$ such that $d^k_{A^k} = 0$
\item[]$\,\,\,8$\hspace*{0.9truecm} If $\nabla f(\tilde x^k)^Td^k < 0$ then
\item[]$\,\,\,9$\hspace*{1.7truecm} Compute a stepsize $\alpha^k\in (0,\alpha^k_{\text{max}}]$ by means of a line search
\item[]$10$\hspace*{0.9truecm} Else
\item[]$11$\hspace*{1.7truecm} Set $\alpha^k = 0$
\item[]$12$\hspace*{0.9truecm} End if
\item[]$13$\hspace*{0.9truecm} Set $x^{k+1}=\tilde x^k + \alpha^k d^k$
\item[]$14$\hspace*{0.1truecm} End for
\par\vspace*{0.1cm}
}
\end{algorithmic}
\end{algorithm}

\subsection{Global Convergence Analysis}\label{sec:glob_conv_simplex}
In this section, we show the global convergence of \ASSIMPL\ to stationary points.
First, we need to describe how we compute the search direction $d^k$ and the stepsize $\alpha^k$ (lines~7 and 9 in Algorithm~\ref{alg:AS_SIMPL}, respectively).

From now on, given any feasible point $x^k$ generated by the algorithm and a feasible direction $d^k$,
we call $\alpha^k_{\text{max}}$ the maximum feasible stepsize that can be used along this direction.
Taking inspiration from~\cite{bertsekas1999nonlinear}, in our framework we require the search direction to be \textit{active-set gradient related}, according to the
following definition:
\begin{definition}\label{def:dirgradrel}
Given an infinite sequence of points $\{x^k\}$ produced by \ASSIMPL,
we say that the sequence of directions $\{d^k\}$ is \textit{active-set gradient related} if,
for any subsequence $\{x^k\}_K$  such that $N(x^k) = \hat N$ for all $k \in K$ and
$\displaystyle{\lim_{k \to \infty, \, k \in K} x^k = x^*}$, where $x^*$
is non-stationary in $\Delta_{\hat N}$, we have that
\begin{align}
& \{d^k_{\hat N}\}_K \text{ is bounded}, \label{dir_grad_rel1} \\[1.4ex]
& \limsup_{k \to \infty, \, k\in K}  \nabla_{\hat N} f(\tilde x^k)^T d_{\hat N}^k < 0, \label{dir_grad_rel2} \\[1.0ex]
& \liminf_{k \to \infty, \, k\in K} \alpha^k_{\text{max}} > 0. \label{dir_grad_rel3}
\end{align}
\end{definition}

By a slight abuse of terminology, in the following we say that a given direction $d^k$ is active-set gradient related if
it is an element of a sequence $\{d^k\}$ satisfying Definition~\ref{def:dirgradrel}.
Moreover, in Subsection~\ref{subsec:dir_simplex}, we will provide some examples of active-set gradient related directions.

For what concerns the computation of the stepsize, a possibility is that of considering the classical Armijo line search
(see, e.g.,~\cite{bertsekas1999nonlinear} and references therein).
This method, which basically performs a successive stepsize reduction, allows to avoid the often considerable
computation associated with an exact line search. Indeed, when dealing with some non-convex problems,
even finding  an approximate local minimizer along the search direction generally requires too many evaluations of the objective
function and possibly the gradient.

The detailed scheme of the Armijo line search is reported in Algorithm~\ref{alg:ls_as}.

\begin{algorithm}[h!]
\caption{\texttt{Armijo line search }} 
\label{alg:ls_as}
\begin{algorithmic}
{\small
\par\vspace*{0.1cm}
\item[]$0$\hspace*{0.5truecm} Choose $\delta \in (0,1)$, $\gamma\in (0,1)$
\item[]$1$\hspace*{0.5truecm} Set initial stepsize $\alpha=\alpha^k_{\text{max}}$
\item[]$2$\hspace*{0.5truecm} While $f(\tilde x^k+ \alpha d^k )> f(\tilde x^k)+\gamma\, \alpha\, \nabla f(\tilde x^k)^T d^k$
\item[]$3$\hspace*{2.0truecm} Set $\alpha = \delta \alpha$
\item[]$4$\hspace*{0.5truecm} End while
\par\vspace*{0.1cm}
}
\end{algorithmic}
\end{algorithm}

Now, we are ready to show that \ASSIMPL\ globally converges to stationary points when
the sequence of directions is active-set gradient related and the Armijo line search is used to compute the stepsize.

\begin{theorem}\label{th:convres}
Let Assumption~\ref{ass:eps} hold. Let $\{x^k\}$ be the sequence of points produced by~\ASSIMPL, where
the sequence of directions $\{d^k\}$ is active-set gradient related and the stepsize $\alpha^k$ is computed using the Armijo line search.
Then, either an integer $\bar k \ge 0$ exists such that $x^{\bar k}$ is a stationary point for problem~\eqref{simplex_prob},
or the sequence $\{x^k\}$ is infinite and every limit point $x^*$ of the sequence is a stationary point for problem~\eqref{simplex_prob}.
\end{theorem}

\begin{proof}
Let $\{x^k\}$ be the sequence produced by~\ASSIMPL\ and let us assume that a stationary point is not produced in a finite number of iterations.
First note that, from the instructions of the algorithm and Proposition~\ref{prop_as}, we can write
\begin{equation*}
f(x^{k+1})\le f(\tilde x^k)\le f(x^k) - CL \norm{\tilde x^k-x^k}^2.
\end{equation*}
Therefore, from the continuity of the objective function and the compactness of the feasible set, it follows that
\begin{gather}
\lim_{k\rightarrow \infty} [f(x^{k+1}) - f(x^k)]  = 0, \label{simpl_f_conv} \\
\lim_{k \to \infty} \|\tilde x^k-x^k\| = 0. \label{simpl_x_tilde_to_x}
\end{gather}

Since the feasible set is compact, then the sequence $\{x^k\}$ attains a limit point $x^*$ and, recalling~\eqref{simpl_x_tilde_to_x},
there exists $K\subseteq \mathbb N$ such that
\begin{equation}\label{convxk}
\lim_{k\rightarrow \infty, \, k\in K} x^k = \lim_{k\rightarrow \infty, \, k\in K} \tilde x^k =  x^*.
\end{equation}
Taking into account the structure of the feasible set, we can characterize a stationary point $x$ using the following condition (see Remark~\ref{remark:kkt}):
\begin{equation*}
\nabla f(x)^T (e_i - x) \ge 0, \quad \forall \ i \in \{1,\dots,n\}.
\end{equation*}
Let $\Phi_i(x)$ be the continuous function defined as
\begin{equation*}
\Phi_i(x) = \max\{0, -\nabla f(x)^T (e_i - x) \}, \quad i = 1,\ldots,n,
\end{equation*}
that measures the violation of the stationarity conditions for a variable $x_i$, $i = 1,\ldots,n$.

By contradiction, we assume that $x^*$ is non-stationary, so that an index $\nu \in \{1,\ldots,n\}$ exists such that
\begin{equation}\label{contr}
\Phi_\nu(x^*) > 0.
\end{equation}
Taking into account that the number of possible different choices of $A^k$ and $N^k$ is
finite, we can find a subset of iteration indices $\bar K \subseteq K$  such that $A^k=\hat A$ and $N^k=\hat N$ for all $k\in \bar K$.

First, suppose that $\nu \in \hat A$. Then, by Definition~\ref{def:activesetest}, we can write
\begin{equation*}
0 \le x^k_\nu \le \epsilon \nabla f(x^k)^T(e_\nu - x^k),
\end{equation*}
so that $\Phi_\nu(x^k) = \max\{0, -\nabla f(x^k)^T (e_\nu - x^k) \} = 0$, for all $k \in \bar K$.
Therefore, from~\eqref{convxk} and the continuity of the function $\Phi_i(\cdot)$, we get a contradiction with~\eqref{contr}.

Then, $\nu$ necessarily belongs to $\hat N$, that is, $x^*_{\hat N}$ is non-stationary in $\Delta_{\hat N}$, where $\Delta_{\hat N}$
is given as in Definition~\ref{def:dirgradrel}.
From that definition and the fact that $d^k_{\hat A} = 0$, we also have that $\{d^k\}_K$ is bounded.
Then, there exists an infinite subset of $K$ (that we still denote by $K$ without loss of generality) such that
\begin{equation}\label{lim_d}
\lim_{k \to \infty \, k \in K} d^k = \bar d.
\end{equation}
Moreover, exploiting again Definition~\ref{def:dirgradrel}, we have that $\eta > 0$ and $M>0$ exist such that
\begin{gather}
\limsup_{k \to \infty, \, k\in K}  \nabla f(\tilde x^k)^T d^k = -\eta, \label{simpl_lim_to_eta} \\
\alpha^k_{\text{max}} \ge M > 0, \quad \forall k \text{ sufficiently large}, \; k \in K \label{alpha_max}.
\end{gather}
From~\eqref{simpl_lim_to_eta}, it follows that $\hat k \in K$ exists such that $\nabla f(\tilde x^k)^T d^k < 0$, for $k \ge \hat k$, $k \in K$.
Then, according to line~9 of Algorithm~\ref{alg:AS_SIMPL}, the Armijo line search computes a value $\alpha^k \in (0,\alpha^k_{\text{max}}]$
for all $k \ge \hat k$, such that
\begin{equation}\label{arm_decr}
f(x^{k+1}) \le f(\tilde x^k) + \gamma \, \alpha^k \, \nabla f(\tilde x^k)^T d^k, \quad \forall k \ge \hat k, \, k \in K,
\end{equation}
or equivalently,
\begin{equation*}
f(\tilde x^k) - f(x^{k+1}) \ge \gamma \, \alpha^k \, \abs{\nabla f(\tilde x^k)^T d^k}, \quad \forall k \ge \hat k, \, k \in K.
\end{equation*}
From~\eqref{simpl_f_conv} and~\eqref{simpl_x_tilde_to_x}, we get that the left-hand side of the above inequality converges to zero for $k \to \infty$,
hence
\begin{equation}\label{simpl_alpha_gd_to_zero}
\lim_{k \to \infty, \, k \in K} \alpha^k\, \abs{\nabla f(\tilde x^k)^T d^k} = 0.
\end{equation}
Using~\eqref{simpl_lim_to_eta}, we obtain that $\displaystyle{\lim_{k \to \infty, \, k \in K} \alpha^k = 0}$.
Taking into account~\eqref{alpha_max}, it follows that there exists $\bar k \in K$, $\bar k \ge \hat k$, such that
\begin{equation*}
\alpha^k < \alpha^k_{\text{max}}, \quad \forall k \ge \bar k, \, k \in K.
\end{equation*}
In other words, for $k \ge \bar k$, $k \in K$, the stepsize $\alpha^k$ cannot be set equal to the maximum feasible stepsize
and, taking into account the line search procedure, we can write
\begin{equation}\label{simpl_alphadelta}
f\Bigl( \tilde x^k + \frac{\alpha^k}{\delta} d^k \Bigr) >
f(\tilde x^k) + \gamma \, \frac{\alpha^k}{\delta}\, \nabla f(\tilde x^k)^T d^k, \quad \forall k \ge \bar k, \, k \in K.
\end{equation}
We can apply the mean value theorem and we have that $\xi_k\in (0,1)$ exists such that
\begin{equation}\label{simpl_meanvaltheo}
f\Bigl( \tilde x^k+ \frac{\alpha^k}{\delta} d^k \Bigr)
= f(\tilde x^k) + \frac{\alpha^k}{\delta}\nabla f\Bigl( \tilde x^k + \xi_k\frac{\alpha^k}{\delta} d^k \Bigr)^T d^k,
\quad \forall k \ge \bar k, \, k \in K.
\end{equation}
By substituting~\eqref{simpl_meanvaltheo} within~\eqref{simpl_alphadelta}, we have
\begin{equation}\label{simpl_meanvaltheo2}
\nabla f\Bigl( \tilde x^k + \xi_k \frac{\alpha^k}{\delta} d^k \Bigr)^T d^k > \gamma\, \nabla f(\tilde x^k)^T d^k,
\quad \forall k \ge \bar k, \, k \in K.
\end{equation}
From~\eqref{convxk}, and exploiting the fact that $\{\xi_k\}$, $\{\alpha^k\}$ and $\{d^k\}$ are bounded, we get
\begin{equation*}
\lim_{k \to \infty, \, k \in K} \tilde x^k + \xi_k \frac{\alpha^k}{\delta} d^k = \lim_{k \to \infty, \, k \in K} \tilde x^k = x^*.
\end{equation*}
Therefore, taking the limits in~\eqref{simpl_alphadelta} and~\eqref{simpl_meanvaltheo}, and taking into account~\eqref{lim_d},
we obtain that
$\nabla f(x^*)^T \bar d \ge \gamma\, \nabla f(x^*)^T \bar d$, or equivalently,
\begin{equation*}
(1-\gamma) \nabla f(x^*)^T \bar d \ge 0.
\end{equation*}
Since $\gamma \in (0,1)$, it follows that $\nabla f(x^*)^T \bar d \ge 0$, contradicting~\eqref{simpl_lim_to_eta}. Hence,
we get $\Phi_i(x^*)=0$, for all $i=1,\dots,n$ and $x^*$ is a stationary point for problem~\eqref{simplex_prob}.
\end{proof}

\begin{remark}
Theorem~\ref{th:convres} holds when using as stepsize in \ASSIMPL\ any value $\alpha^k\in (0,\alpha^k_{\text{max}}]$ such that
\[
f(\tilde x^k + \alpha^k d^k)\le f(\tilde x^k + \alpha_A^k d^k),
\]
where $\alpha_A^k$ is the value computed by the Armijo line search.
It follows from the fact that, if the above relation is satisfied, then~\eqref{arm_decr} holds, as well as all the subsequent
steps in the proof.

In particular, this implies that Theorem~\ref{th:convres} holds under the assumption that
the stepsize is computed in \ASSIMPL\ by means of an exact line search, that is, $\alpha^k$ is computed as
\begin{equation*}
\alpha^k \in \argmin_{\alpha \in (0,\alpha^k_{\text{max}}]} f(\tilde x^k + \alpha d^k).
\end{equation*}
\end{remark}

\subsection{Active-Set Gradient Related Directions in \ASSIMPL}\label{subsec:dir_simplex}
As shown in Theorem~\ref{th:convres}, to ensure global convergence of \ASSIMPL\ to stationary points
we need a sequence of directions that satisfies Definition~\ref{def:dirgradrel}.
Moreover, from the instructions of the algorithm, $d^k_{A^k}$ must be zero at each iteration, so we only need to focus on the computation of $d^k_{N^k}$.

As examples of $d^k_{N^k}$ that can be used in practice, we consider Frank-Wolfe-type and projected gradient directions:
\medskip\par
\begin{description}
\item[(FW)] \textit{Frank-Wolfe direction}:
    \begin{equation}\label{d_fw}
    d^{\text{FW}}_{N^k} = (e_{\hat\imath} - \tilde x^k)_{N^k}, \quad \hat\imath \in \argmin_{i\in N^k} \bigl\{ \nabla_i f(\tilde x^k) \bigr\};
    \end{equation}
\item[(AFW)] \textit{away-step Frank-Wolfe direction}:
    \begin{equation*}
    d^{\text{AFW}}_{N^k} =
    \begin{cases}
    d^{\text{FW}}_{N^k}, \quad & \text{ if } \nabla_{N^k} f(\tilde x^k)^T d^{\text{FW}}_{N^k} \le \nabla_{N^k} f(\tilde x^k)^T d^{\text{A}}_{N^k}, \\
    d^{\text{A}}_{N^k}, \quad & \text{ otherwise},
    \end{cases}
    \end{equation*}
    where
    \begin{equation}\label{d_a}
    d^{\text{A}}_{N^k} = (\tilde x^k - e_{\hat\jmath})_{N^k}, \quad  \hat\jmath \in \argmax_{j\in  N^k_0} \bigl\{ \nabla_j f(\tilde x^k) \bigr\}
    \end{equation}
    and $N^k_0 = \{j\in N^k \colon \tilde x_j^k > 0\}$;
    \smallskip\par
\item[(PG)] \textit{projected gradient direction}:
    \begin{equation*}
    d^{\text{PG}}_{N^k} = \Bigl(P\bigl(\tilde x^k -s\nabla f(\tilde x^k)\bigr)_{\Delta_{N^k}}-\tilde x^k\Bigr)_{N^k},
    \end{equation*}
    where $s>0$ is a fixed scalar.
\end{description}
\medskip\par

In the following, we will refer to $d^{\text{FW}}$, $d^{\text{AFW}}$ and $d^{\text{PG}}$
when the subdirection $d^k_{N^k}$ is chosen according to the Frank-Wolfe (FW), the away-step Frank-Wolfe (AFW)
or the projected gradient (PG) rule, respectively.
We now show that the three considered directions satisfy Definition~\ref{def:dirgradrel}.
\begin{proposition}\label{prop:grad_dir}
Given one rule among (FW), (AFW) and (PG) for the computation of $d^k_{N^k}$,
the resulting sequence of directions $\{d^k\}$ generated by \ASSIMPL\ is active-set gradient related.
\end{proposition}

\begin{proof}
Since $\Delta$ is compact, it is easy to see that all the considered directions are bounded.
Considering some of the ideas reported in~\cite{bertsekas1999nonlinear} (see chapter 2)
and the properties related to the active-set estimate, we now prove that those directions
also satisfy~\eqref{dir_grad_rel2} and~\eqref{dir_grad_rel3}.
Let $\{x^k\}_K$ be a subsequence such that $N^k = \hat N$ for all $k \in K$
and $\displaystyle{\lim_{k \to \infty, \, k \in K} x^k_{\hat N} = x^*_{\hat N}}$,
where $x^*$ is non-stationary in $\Delta_{\hat N}$.
We consider the different cases:
\begin{description}
\item[(FW)] By definition of the index $\hat \imath$ given in~\eqref{d_fw}, it is easy to see that, for all $k \ge 0$, we have
    $\nabla_{\hat \imath} f(\tilde x^k) \le \nabla_{N^k} f(\tilde x^k)^T x_{N^k}$, $\forall x \in \Delta_{N^k}$.
    Thus,
    \[
    \begin{split}
    \nabla_{\hat N} f(\tilde x^k)^T d_{\hat N}^k & = \nabla_{\hat N} f(\tilde x^k)^T (e_{\hat \imath}-\tilde x^k)_{\hat N} \\
                                                 & \le \nabla_{\hat N} f(\tilde x^k)^T (x-\tilde x^k)_{\hat N} = \nabla f(\tilde x^k)^T (x-\tilde x^k),
                                                 \quad \forall x \in \Delta_{\hat N},
    \end{split}
    \]
    where the last equality follows from the fact that $x_{A^k} = 0$ for all $x \in \Delta_{N^k}$ and $\tilde x^k_{A^k} = 0$.
    Passing to the limit, we obtain
    \begin{equation*}
    \limsup_{k \to \infty, \, k\in K}  \nabla_{\hat N} f(\tilde x^k)^T d_{\hat N}^k
    \le \nabla f(x^*)^T (x-x^*), \quad \forall x \in \Delta_{\hat N}.
    \end{equation*}
    Since $x^*$ is non-stationary in $\Delta_{\hat N}$, we have that
    \[
    \min_{x \in \Delta_{\hat N}} \nabla f(x^*)^T (x-x^*) < 0.
    \]
    Therefore, combining the above two inequalities,~\eqref{dir_grad_rel2} holds.
    Since $\alpha^k_{\text{max}} = 1$ at every iteration, also~\eqref{dir_grad_rel3} is trivially satisfied.
\item[(AFW)] Since, by definition, $\nabla f(\tilde x^k)^Td^{\text{AFW}} \le \nabla f(\tilde x^k)^Td^{\text{FW}}$,
    we can repeat the same reasoning given above for the (FW) case and~\eqref{dir_grad_rel2} holds.
    To prove~\eqref{dir_grad_rel3}, by contradiction let us assume that an infinite subset of $K$
    (that we still denote by $K$ for the sake of simplicity) exists such that
    \begin{equation}\label{limalpha}
    \lim_{k \to \infty, \, k\in K} \alpha_{\text{max}}^k =0.
    \end{equation}
    Recalling the definition of $d^{\text{AFW}}$, the case we need to analyze is the one where we get an infinite subsequence of
    away-step directions in $\hat N$ (as $\alpha^k_{\text{max}} = 1$ for Frank-Wolfe directions).
    So, we assume that an infinite subset $\tilde K \subseteq K$ exists such that
    \begin{equation*}
    d^k_{\hat N} = d^{\text{A}}_{\hat N}, \quad  \forall k \in \tilde K.
    \end{equation*}
    We have that $\alpha^k_{\text{max}}= \displaystyle\frac{\tilde x_{\hat\jmath}^k}{1 - \tilde x_{\hat\jmath}^k}$, for all $k \in \tilde K$,
    where $\hat\jmath $ is the index computed according to~\eqref{d_a}.
    Since the number of indices in $\hat N$ is finite, we can consider an infinite subset of $\tilde K$
    (that we still denote by $\tilde K$ for the sake of simplicity) where the index $\hat\jmath$ is fixed.
    Taking into account~\eqref{limalpha}, it is easy to see that $\{\tilde x^k_{\hat\jmath}\}_{\tilde K} \to 0$.
    Using~\eqref{simpl_x_tilde_to_x}, we get
    \begin{equation}\label{limxkj}
    \lim_{k \to \infty, \, k\in \tilde K} x^k_{\hat\jmath} = 0.
    \end{equation}
    Moreover, from~\eqref{dir_grad_rel2}, \eqref{simpl_x_tilde_to_x} and the continuity of $\nabla f(x)$, we can write
    \[
    \limsup_{k \to \infty, \, k\in K} \nabla_{\hat N} f(x^k)^T d^k_{\hat N} = -\eta.
    \]
    Exploiting the fact that $\nabla f(x^k)^T d^k = \nabla_{\hat N} f(x^k)^T d^k_{\hat N}$ by definition of $d^k$, we obtain
    \begin{equation}\label{limdirder}
    \begin{split}
    -\eta & = \limsup_{k \to \infty, \, k\in K} \nabla f(x^k)^T d^k \\[1ex]
          & = \limsup_{k \to \infty, \, k\in K} [\nabla f(x^k)^T (\tilde x^k - x^k) + \nabla f(x^k)^T (x^k - e_{\hat \jmath})] \\[1ex]
          & = \limsup_{k \to \infty, \, k\in K} \nabla f(x^k)^T (x^k - e_{\hat \jmath}),
    \end{split}
    \end{equation}
    where the last equality follows again from~\eqref{simpl_x_tilde_to_x}.
    From~\eqref{limxkj} and~\eqref{limdirder}, an index $\tilde k \in \tilde K$ exists such that, for all $k \ge \tilde k$, $k \in \tilde K$, we have
    \begin{gather*}
    \nabla f(x^k)^T (x^k - e_{\hat \jmath}) \le -\frac{\eta} 2, \\
    x^k_{\hat\jmath} \le \epsilon \frac \eta 2.
    \end{gather*}
    Therefore,
    \begin{equation*}
    x^k_{\hat\jmath} \le \eps \nabla f(x^k)^T (e_{\hat\jmath} - x^k), \quad \forall k \ge \tilde k, \, k \in \tilde K.
    \end{equation*}
    Recalling~\eqref{as1}, this implies that $\hat \jmath \notin \hat N$ and, considering the definition of $\hat \jmath$ in~\eqref{d_a}, we get a contradiction.
\item[(PG)] Let us define $\hat x^k = P\bigl(\tilde  x^k-s\nabla f(\tilde x^k)\bigr)_{\Delta_{N^k}}$, so that $d^k = \hat x^k-\tilde x^k$.
    By continuity of the projection operator, we get
    \[
    \lim_{k \to \infty, \, k\in K} \hat x^k = P\bigl( x^* - s\nabla f( x^*)\bigr)_{\Delta_{\hat N}}.
    \]
    From the properties of the projection, we have
    \[
    (\tilde x^k - s\nabla f(\tilde x^k)-\hat x^k)^T(x -\hat x^k) \le 0, \quad \forall x\in \Delta_{N^k}.
    \]
    If we choose $x=\tilde x^k$ in the above inequality, we can write
    \begin{equation}\label{pg_ineq}
    \nabla f(\tilde x^k)^T d^k = \nabla f(\tilde x^k)^T (\hat x - \tilde x^k) \le - \frac 1 s \|\hat x^k - \tilde x^k\|^2
    = - \frac 1 s \|d^k\|^2, \quad \forall k \ge 0.
    \end{equation}
    Since $d^k_{A^k} = 0$ for all $k \ge 0$, taking the limit we have
    \[
    \begin{split}
    \limsup_{k \to \infty, \, k\in K} \nabla_{\hat N} f(\tilde x^k)^T d_{\hat N}^k
    & = \limsup_{k \to \infty, \, k\in K} \nabla f(\tilde x^k)^T d^k \\
    & \le -\frac 1 s \bigl\|P\bigl(x^* -s\nabla f(x^*)\bigr)_{\Delta_{\hat N}}- x^*\bigr\|^2.
    \end{split}
    \]
    From the fact that $x^*$ is non-stationary in $\Delta_{\hat N}$, it follows that
    \[\bigl\|P\bigl(x^* -s\nabla f(x^*)\bigr)_{\Delta_{\hat N}}- x^*\bigr\| > 0.\] Therefore,
    \[
    \limsup_{k \to \infty, \, k\in K} \nabla_{\hat N} f(\tilde x^k)^T d_{\hat N}^k < 0,
    \]
    implying that~\eqref{dir_grad_rel2} holds.
    Finally, since $\alpha^k_{\text{max}} = 1$ at every iteration, also~\eqref{dir_grad_rel3} is trivially satisfied.
\end{description}
 \end{proof}

\begin{remark}\label{remark:dir}
Since we set $\tilde x^k_{A^k} = 0$ at any iteration $k$, it is straightforward to verify that,
when $d^k$ is computed according to (FW), (AFW) or (PG) rule, $\nabla f(\tilde x^k)^T d^k < 0$
if and only if $\tilde x^k$ is non-stationary on $\Delta_{N^k}$.
Equivalently, $\nabla_{N^k} f(\tilde x^k)^T d^k_{N^k} < 0$ if and only if $\tilde x^k_{N^k}$ is non-stationary on the subspace variable $N^k$.
\end{remark}

\section{Convergence Rate Analysis}\label{sec:convrate}
In this section, we analyze the convergence rate of \ASSIMPL\ when one rule among (FW), (AFW) and (PG) is used
for the computation of $d^k$. More specifically, we focus on particular classes of non-convex problems
(i.e., problems satisfying some specific assumptions we will make later on),
and report linear convergence results for our framework when using those three directions.
The results are asymptotic since they exploit the properties of the active-set estimate given in Theorem~\ref{th:estim},
and these properties hold only in a neighborhood of stationary points.
Summarizing, on the one hand, we get asymptotic linear rate, but, on the other hand, our results hold for non-convex objective functions.

We make an assumption that is pretty common when analyzing the convergence
rate of algorithms (see, e.g.,~\cite{ortega2000iterative}), and quite reasonable, taking into account the results
reported in the previous section.

\begin{assumption}\label{ass:rate}
Let $\{x^k\}$ be the infinite sequence generated by \ASSIMPL. We have that
\begin{equation*}
\lim_{k \to \infty} x^k=x^*,
\end{equation*}
where $x^*$ is a stationary point of problem~\eqref{simplex_prob}.
\end{assumption}

From now on, we denote with $\bar I$ the set $\{1,\ldots,n\}$.
We also denote  with $\bar A$ and $\bar N$ the index sets defined in~\eqref{AS} and~\eqref{NAS}, respectively, and with
\[
N^+:=\bar N \cup \{i\in \bar I \colon x^*_i =0, \, \mu^*_i =0\} \quad \text{and} \quad A^+:= \bar I \setminus N^+= \{i\in \bar I \colon x^*_i =0, \, \mu^*_i >0\}.
\]

\subsection{Linear Convergence of Active-Set Frank-Wolfe}
Here we show that, when the Frank-Wolfe direction (FW) is embedded in our active-set framework, one
can get asymptotic linear convergence without making the classic assumptions (see, e.g.,~\cite{guelat:1986})
needed for proving linear convergence rate of the classical Frank-Wolfe method, that is:
\begin{itemize}
 \item optimal solution in the interior of the feasible set,
 \item strongly convex objective function.
\end{itemize}
As we will see, those assumptions are replaced by strict complementarity in the optimal
solution and strong convexity on $\Delta_{\bar N}$, respectively. Again,
we remark that the results are asymptotic, but they do not require convexity assumptions of the objective function on the whole $\Delta$.
Moreover, we obtain pretty tight convergence rate constants, that get much tighter than those obtained with the classical Frank-Wolfe method
as the final solution is sparse.
So, we can consider the result as a good trade-off in the end.

Before reporting the theoretical results related to the active-set Frank-Wolfe (i.e., \ASSIMPL~with $d^k$ computed according to (FW) rule),
we need to introduce some constants, which follow from those used in~\cite{lacoste:2013}, adapted to our purposes.
Given an index subset $I \subseteq \bar I$, we define:
\begin{gather*}
C_f(I) := \sup_{\substack{x,s \in \Delta_I, \\ \alpha \in (0,1], \\ y = x + \alpha (s-x)}}
     \frac 2 {\alpha^2} \Bigl[ f(y) - f(x) - \nabla f(x)^T (y-x) \Bigr], \\
\mu_f(I) := \inf_{\substack{x \in \Delta_I \setminus \{x^*\}, \\ \alpha \in (0,1], \\ \bar s = \bar s(x,x^*,\Delta), \\ y = x + \alpha(\bar s-x)}}
    \frac 2 {\alpha^2}\Bigl[ f(y) - f(x) - \nabla f(x)^T (y-x) \Bigr],
\end{gather*}
where $\bar s(x,x^*,\Delta) := \text{ray}(x,x^*) \cap \partial (\Delta)$ and $\text{ray}(x,x^*)$ is the ray from $x$ to $x^*$.
The curvature constant $C_f(I)$, which  measures the non-linearity of the objective function in the subspace $\Delta_I$,
is needed to give a quadratic upper bound on the objective function.
The strong convexity constant $\mu_f(I)$, which measures the strong convexity of the objective function
on $\Delta_I$ (and can be interpreted as the lower curvature of the function), is used to give a quadratic lower bound instead
(see~\cite{lacoste:2013} for further details).

\begin{remark}
The main difference between the constants given above and those introduced in~\cite{lacoste:2013}
is that ours are restricted to a particular subspace.
Moreover, for any index subset $I \subseteq \bar I$, it is easy to see that
\begin{equation}\label{ineq_rate_1}
\mu_f(\bar I) \le \mu_f(I) \le C_f(I) \le C_f(\bar I).
\end{equation}
\end{remark}

Now, we are ready to state linear convergence rate of \ASSIMPL\ when (FW) rule is used to compute the search direction.

\begin{theorem}\label{th:linrateasfw}
Let Assumption~\ref{ass:eps} and~\ref{ass:rate} hold, let $f(x)$ be strongly convex on $\Delta_{\bar N}$,
and let us assume that strict complementarity holds at $x^*$.
Let us further assume that $d^k$ is computed by (FW) rule and that the exact line search is used.

Then, there exists $\bar k$ such that
\begin{equation*}
f(x^{k+1}) - f(x^*) \le \bigl(1 - \rho^{AS-FW}\bigr) \bigl[ f(x^k) - f(x^*) \bigr], \quad \forall k \ge \bar k,
\end{equation*}
where
\begin{equation*}
\rho^{\text{AS-FW}} = \min \biggl\{ \frac 1 2, \frac{\mu_f(\bar N)}{C_f (\bar N)} \biggr\}.
\end{equation*}
\end{theorem}
\begin{proof}
From Theorem~\ref{th:estim}, exploiting the fact that strict complementarity holds at $x^*$, for sufficiently large $k$ we have that
\begin{equation*}
N(x^k) = N(\tilde x^k) = \bar N \quad \text{ and } \quad A(x^k) = A(\tilde x^k) = \bar A.
\end{equation*}
From the instructions of \ASSIMPL, for sufficiently large $k$ we have $\tilde x^k_{\bar A} = x^k_{\bar A} = 0$,
implying that $\tilde x^k = x^k$. Then, for sufficiently large $k$ the minimization is restricted to the variable subspace $N^k = \bar N$.
Since the search direction $d^k$ is computed according to (FW) rule, the rest of the proof follows
by repeating the same arguments of the proof given for Theorem~$3$ in~\cite{lacoste:2013},
observing that $\mu_f(\bar N) > 0$ and $C_f(\bar N) < \infty$ under the hypothesis we made.
\end{proof}

\begin{remark}
From~\eqref{ineq_rate_1}, it follows that the smaller $\bar N$ (i.e., the sparser $x^*$),
the better the convergence rate of \ASSIMPL.
Moreover,
\begin{equation*}
\rho^{\text{AS-FW}} \ge \min \biggl\{ \frac 1 2, \frac{\mu_f(\bar I)}{C_f (\bar I)} \biggr\} = \rho^{\text{FW}},
\end{equation*}
where $\rho^{\text{FW}}$ is the constant given in~\cite{lacoste:2013} for the convergence rate of the standard Frank-Wolfe method.
\end{remark}

\subsection{Linear Convergence of Active-Set Away-Step Frank-Wolfe}
In this subsection, we prove that active-set away-step Frank-Wolfe
(i.e., \ASSIMPL~with $d^k$ computed according to (AFW) rule) asymptotically
converges at linear rate. We can prove the result without making the strong convexity assumption (see, e.g.,~\cite{guelat:1986})
needed for proving linear convergence rate of the classical away-step Frank-Wolfe method.
As we will see, that assumption is replaced by strong convexity of the objective function on $\Delta_{N^+}$.
Similarly to the (FW) direction, here we get asymptotic results,
but we do not need strong convexity assumptions of the objective function on the whole $\Delta$
and we obtain pretty tight convergence rate constants that depend on the sparsity of the final solution.
Again, we can consider the result as a good trade-off in the end.

Given an index subset $I \subseteq \bar I$, we define the following two constants,
which follow from those used in~\cite{lacoste:2015}, adapted to our purposes:
\begin{gather*}
C^{\Delta}_f(I) := \sup_{\substack{x,s,v \in \Delta_I \\ \alpha \in (0,1], \\ y = x + \alpha (s-v)}}
     \frac 2 {\alpha^2} \Bigl[ f(y) - f(x) - \alpha\nabla f(x)^T (s-v) \Bigr], \label{Cf}\\
\mu^{\Delta}_f(I) := \inf_{x \in \Delta_I } \; \inf_{\substack{\hat x \in \Delta_I \\ \nabla f(x)^T (\hat x - x) < 0}}
     \frac 2 {{\alpha^{\Delta}_I(x,\hat x)}^2}\Bigl[ f(\hat x) - f(x) - \nabla f(x)^T (\hat x-x) \Bigr],
\end{gather*}
where
\begin{gather*}
\alpha^{\Delta}_I(x,\hat x) := \frac{\nabla f(x)^T (\hat x-x)}{\nabla f(x)^T (s_I(x)-v_I(x))}, \\
s_I(x) := e_{\hat \imath}, \quad \hat \imath \in \argmin_{i \in I} \{ \nabla_i f(x) \}, \\
v_I(x) := e_{\hat \jmath}, \quad \hat \jmath \in \argmax_{j \in I \colon x_j > 0} \{ \nabla_j f(x) \}.
\end{gather*}
These two new constants are motivated in the analysis by the fact that both Frank-Wolfe and away-step directions
are used (see~\cite{lacoste:2015} for further details).

\begin{remark}
Also in this case, the main difference between the constants given above and those introduced in~\cite{lacoste:2015}
is that ours are restricted to a particular subspace.
Moreover, for any index subset $I \subseteq \bar I$, it is easy to see that the following inequalities hold:
\begin{equation}\label{ineq_rate_2}
\mu^{\Delta}_f(\bar I) \le \mu^{\Delta}_f(I) \le C^{\Delta}_f(I) \le C^{\Delta}_f(\bar I).
\end{equation}
\end{remark}

Theorem~$8$ in~\cite{lacoste:2015} shows, for the standard away-step Frank-Wolfe method,
that the quantity $[f(x^k)-f(x^*)]$ decreases linearly
at each iteration $k$ that is not a so-called \textit{drop step}.
Iteration $k$ is a drop step when the stepsize $\alpha^k = \alpha^k_{\text{max}}<1$ and the number of zero components in $x^{k+1}$
increases by one.
In the convergence rate analysis, these iterations are troublesome since a geometric decrease of $[f(x^k)-f(x^*)]$ cannot be guaranteed.

In our context, these definitions apply when considering the computation of $x^{k+1}$ from $\tilde x^k$ and, as to be shown in the next theorem,
we can still guarantee that the quantity $[f(x^k)-f(x^*)]$ decreases linearly at each iteration $k$ that is a \textit{good step} (i.e., not a drop step)
with tighter constants (that depend on the sparsity of the optimal solution).

\begin{theorem}\label{th:linrateasvar}
Let Assumption~\ref{ass:eps} and~\ref{ass:rate} hold, let $f(x)$ be strongly convex on $\Delta_{N^+}$,
with $\nabla f(x)$ Lipschitz continuous on $\Delta_{N^+} + (\Delta_{N^+} - \Delta_{N^+})$ (in the Minkowski sense).
Let us further assume that $d^k$ is computed by (AFW) rule and that the exact line search is used.

Then, there exists $\bar k$ such that, for every iteration $k \ge \bar k$ that is a good step (i.e., it is not a drop step), we have
\begin{equation*}
f(x^{k+1}) - f(x^*) \le (1 - \rho^{\text{AFW}}) \bigl[ f(x^k) - f(x^*) \bigr],
\end{equation*}
where
\begin{equation}\label{rho}
\rho^{\text{AFW}} = \dfrac{\mu^{\Delta}_f(N^+)}{4C_f^{\Delta} (N^+)},
\end{equation}
Moreover, for $k \ge \bar k$, we have that at most $\abs{N^+}-1$ drop steps can be performed in between two good steps.
\end{theorem}

\begin{proof}
First, we observe that Theorem~\ref{th:estim} implies that
an iteration $\tilde k$ exists such that $A^k \supseteq A^+$ and $N^k \subseteq N^+$ for $k \ge \tilde k$.
Now, we show that there exists $\bar k \ge \tilde k$ such that
\begin{enumerate}[(i)]
\item $x^k_{A^+} = \tilde x^k_{A^+} = 0$, for all $k \ge \bar k$;
\item $\nabla f(\tilde x^k)^T d^k < 0$, for all $k \ge \bar k$;
\item $\displaystyle{x^* \in \argmin_{x\in \Delta_{N^k}}~ f(x)}$, for all $k \ge \bar k$.
\end{enumerate}
Point~(i) follows from the instructions of the algorithm and the fact that $A^k \supseteq A^+$, for $k \ge \tilde k$.
To prove point~(ii), we proceed by contradiction.
We assume that an infinite subsequence $\{\tilde x^k\}_K$ exists such that $\nabla f(\tilde x^k)^T d^k = 0$ for all $k \in K$.
Recalling Remark~\ref{remark:dir}, this means that
$\tilde x^k$ is stationary over $\Delta_{N^k}$ (but $\tilde x^k$ is not stationary over $\Delta$), for all $k \in K$.
Since $N^k \subseteq N^+$ for $k \ge \tilde k$ and $f(x)$ is strongly convex on $\Delta_{N^+}$, there exists a unique point satisfying stationarity
conditions over $\Delta_{N^k}$ for $k \ge \tilde k$. Taking into account that $A^k$ and $N^k$ are subsets of a finite set of indices and $\tilde x_{A^k}^k=0$, we
have that, after a finite number of iterations, the algorithm should cycle. This cannot be possible as we guarantee a strict decrease
in the objective function at each iteration.
Point~(iii) follows from the fact that $A^k \subseteq \bar A$ for all $k\ge \tilde k$ and $f(x)$ is strongly convex on $\Delta_{N^+}$.

Consequently, recalling that $d^k_{A^k} = 0$, for $k \ge \bar k$ the minimization is restricted to the variable subspace $N^k \subseteq N^+$.
We can thus repeat the same arguments of the proof given for Theorem~$8$ in~\cite{lacoste:2015} to provide the following bound for all $k \ge \bar k$:
\begin{equation*}
f(x^{k+1}) - f(x^*) \le (1 - \rho^{\text{AFW}}) \bigl[ f(\tilde x^k) - f(x^*) \bigr]
                    \le (1 - \rho^{\text{AFW}}) \bigl[ f(x^k) - f(x^*) \bigr],
\end{equation*}
where the last inequality follows from the fact that $f(\tilde x^k) \le f(x^k)$.
Moreover, we have that $\mu^{\Delta}_f(N^+) > 0$ and $C^{\Delta}_f(N^+) < \infty$ under the hypothesis we made.

Finally, to bound the number of iterations for which $k$ is not a good step,
we need to consider those iterations such that $\alpha^k = \alpha^k_{\text{max}} < 1$, for $k \ge \bar k$.
The fact that $\alpha^k_{\text{max}} < 1$ implies that $d^k = d^{\text{A}}$.
Consequently, when $\alpha^k = \alpha^k_{\text{max}}$, we have that $x^{k+1}_{\hat \jmath} = 0$, where $\hat \jmath$ is the index computed
according to~\eqref{d_a}. In other words, the number of zero components in $x^{k+1}$ increases by $1$
(since $d^k_i = 0$ for all $i$ such that $\tilde x^k_i = 0$, i.e., the away-step direction does not change zero components).
From the instructions of the algorithm, we also have that the number of zero components in $\tilde x^{k+1}$ cannot decrease from $x^{k+1}$.
Combining these observations with the fact that $\tilde x^k_{A^+} = 0$ for all $k \ge \bar k$,
we conclude that after at most $\abs{N^+}-1$ iterations with $\alpha^k = \alpha^k_{\text{max}} < 1$,
a point $\tilde x^k$ with $n-1$ zero components is produced. Of course, we cannot further increase the number of zero components. \end{proof}

\begin{remark}
From~\eqref{ineq_rate_2}, it follows that the smaller $N^+$, the better the convergence rate of \ASSIMPL.
Moreover,
\begin{equation*}
\rho^{\text{AS-AFW}} \ge \frac{\mu^{\Delta}_f(\bar I)}{4C_f^{\Delta} (\bar I)} = \rho^{\text{AFW}},
\end{equation*}
where $\rho^{\text{AFW}}$ is the constant given in~\cite{lacoste:2015} for the convergence rate of
the standard away-step Frank-Wolfe.
Furthermore, also the upper bound on the number of bad steps between two good steps
depends on the cardinality of $N^+$ (for sufficiently large $k$).
We would like to recall that, in the standard away-step Frank-Wolfe, this value is equal to $n-1$.
\end{remark}

\subsection{Linear Convergence of Active-Set Projected Gradient}
In this subsection, we prove that the active-set Projected Gradient (i.e., \ASSIMPL~with $d^k$
computed according to (PG) rule) asymptotically converges at a linear rate,
under the assumption that the objective function is strictly convex on $\Delta_{N^+}$.
We follow the same arguments of the proof of Theorem~$3.1$ in~\cite{LuoTseng:1992}.
First, we need to give two additional results, stated in  Lemma~\ref{lemma:pg1} and Lemma~\ref{lemma:pg2}.

\begin{lemma}\label{lemma:pg1}
Let $\{x^k\}$ be the sequence produced by Algorithm~\ref{alg:AS_SIMPL}.
Then, there exists $\bar k$ such that
\[
f(\tilde x^k)- f(x^*) \le \dfrac L 2 \|\tilde x^k - x^*\|^2, \quad \forall k \ge \bar k.
\]
\end{lemma}
\begin{proof}
From the Lipschitz continuity of the gradient, for all $k \ge 0$ we can write
\[
f(\tilde x^k)- f(x^*) \le \nabla f(x^*)^T(\tilde x^k - x^*) + \frac{L}{2}\|\tilde x^k - x^*\|^2.
\]
From Theorem~\ref{th:estim}, an iteration $\bar k$ exists such that $A^k \supseteq A^+$ and $\bar N \subseteq N^k$ for all $k \ge \bar k$.
Hence, from the first-order necessary optimality conditions, $\nabla_i f(x^*) = \lambda^*$ for all $i\in N^k$ and for all $k\ge \bar k$.
Since $\tilde x_i^k = x_i^* = 0$ for all $i \in A^k$ and for all $k\ge \bar k$, we get
\[
\nabla f(x^*)^T(\tilde x^k - x^*) = \sum_{i\in N^k} \lambda^* (\tilde x^k - x^*)_i = 0, \quad \forall k \ge \bar k,
\]
where the last equality follows from the feasibility of $\tilde x^k$ and $x^*$. Therefore, for all $k \ge \bar k$ we obtain
$f(\tilde x^k)- f(x^*) \le L/2 \|\tilde x^k - x^*\|^2$.
\end{proof}

\begin{lemma}\label{lemma:pg2}
Let $\{x^k\}$ be the sequence produced by Algorithm~\ref{alg:AS_SIMPL},
where $d^k$ is computed by (PG) rule and the Armijo line search is used.
Then, at any iteration $k$ such that $\nabla f(\tilde x^k)^T d^k < 0$ we have
\begin{gather}
\alpha^k \ge R, \label{alphabounded} \\[1ex]
f(\tilde x^k) - f(x^{k+1}) \ge \frac{\gamma R \min\{1,s\}^2}s \bigl\|P\bigl(\tilde x^k - \nabla f(\tilde x^k)\bigr)_{\Delta_{N^k}} - \tilde x^k\bigr\|^2, \label{f_decr}
\end{gather}
where $\displaystyle{R= \min\Bigl\{1,\frac{2 \delta (1-\gamma)}{s L}\Bigr\}}$ and $\delta$, $\gamma$ are the parameters used in the Armijo line search.
\end{lemma}
\begin{proof}
Let $k$ be an iteration such that $\nabla f(\tilde x^k)^T d^k < 0$.
Repeating the same reasonings done in Proposition~\ref{prop:grad_dir}, we obtain~\eqref{pg_ineq}, i.e.,
\begin{equation}\label{pg_ineq2}
\nabla f(\tilde x^k)^T d^k \le - \frac 1 s \|d^k\|^2.
\end{equation}

First, we prove~\eqref{alphabounded}.
From the Lipschitz continuity of $\nabla f(x)$, we can write
\begin{equation*}
\begin{split}
f(x^{k+1})- f(\tilde x^k) & \le \nabla f(\tilde x^k)^T (x^{k+1} - \tilde x^k) + \frac L 2 \|x^{k+1} - \tilde x^k \|^2\\
& = \alpha^k \nabla f(\tilde x^k)^T d^k + \frac L 2 (\alpha^k)^2 \|d^k \|^2,
\end{split}
\end{equation*}
Combining the above inequality with~\eqref{pg_ineq2}, we obtain
\[
\begin{split}
f(x^{k+1})- f(\tilde x^k) & \le \alpha^k\nabla f(\tilde x^k)^T d^k - \frac{sL}{2} (\alpha^k)^2 \nabla f(\tilde x^k)^T d^k \\
                          & = \alpha^k \Bigl(1 - \frac{sL\alpha^k}{2}\Bigr) \nabla f(\tilde x^k)^T d^k.
\end{split}
\]
Therefore, recalling that $\alpha$ is multiplied by $\delta$ when is decreased in the line search procedure (see step~3 of Algorithm~\ref{alg:ls_as}),
to satisfy the criterion within the Armijo line search and the feasibility we need to have $\alpha^k \ge R$, and then~\eqref{alphabounded} holds.

Now, we prove~\eqref{f_decr}.
From the Armijo line search,  we have $f(\tilde x^k) - f(x^{k+1}) \ge -\gamma \alpha^k \nabla f(\tilde x^k)^T d^k$. Recalling~\eqref{pg_ineq2}, we obtain
\[
f(\tilde x^k) - f(x^{k+1}) \ge \frac{\gamma \alpha^k}s \|d^k\|^2.
\]
Using the fact that
\[
\norm{d^k} = \bigl\|P\bigl(\tilde x^k - s\nabla f(\tilde x^k)\bigr)_{\Delta_{N^k}} - \tilde x^k\bigr\| \ge
\min\{1,s\} \bigl\|P\bigl(\tilde x^k - \nabla f(\tilde x^k)\bigr)_{\Delta_{N^k}} - \tilde x^k\bigr\|
\]
(see proof of Theorem 4.1 in~\cite{LuoTseng:1992} for the above inequality), we get
\[
f(\tilde x^k) - f(x^{k+1}) \ge \frac{\gamma \alpha^k \min\{1,s\}^2}s \bigl\|P\bigl(\tilde x^k - \nabla f(\tilde x^k)\bigr)_{\Delta_{N^k}} - \tilde x^k\bigr\|^2.
\]
Combining this inequality with~\eqref{alphabounded}, we obtain that~\eqref{f_decr} holds.
\end{proof}

\begin{theorem}\label{th:linrateaspg}
Let $\{x^k\}$ be the sequence produced by Algorithm~\ref{alg:AS_SIMPL},
where $d^k$ is computed by (PG) rule and the Armijo line search is used.
Let Assumption~\ref{ass:eps} and~\ref{ass:rate} hold, and let $f(x)$ be strictly convex on $\Delta_{N^+}$.

Then, there exists $\bar k$ such that
\begin{equation*}
f(x^{k+1}) - f(x^*) \le (1-\rho^{\text{AS-PG}}) \bigl[ f(x^k) - f(x^*) \bigr], \quad \forall k \ge \bar k,
\end{equation*}
with $\rho^{\text{AS-PG}} > 0$.
\end{theorem}
\begin{proof}
From Theorem~\ref{th:estim}, we have that $A^k \supseteq A^+$ and $N^k \subseteq N^+$ for sufficiently large $k$.
Reasoning as in the first part of the proof of Theorem~\ref{th:linrateasvar}, we claim that there exists an iteration $\bar k$ such that
\begin{enumerate}[(i)]
\item $x^k_{A^+} = \tilde x^k_{A^+} = 0$, for all $k \ge \bar k$;
\item $\nabla f(\tilde x^k)^T d^k < 0$, for all $k \ge \bar k$;
\item $\displaystyle{x^* \in \argmin_{x\in \Delta_{N^k}}~ f(x)}$, for all $k \ge \bar k$.
\end{enumerate}
Hence, by Theorem 2.1 in~\cite{LuoTseng:1992}, for sufficiently large $k$ we have
\begin{equation}\label{pg_ineq3}
\|\tilde x^k - x^*\| \le \tau \bigl\|P\bigl(\tilde x^k - \nabla f(\tilde x^k)\bigr)_{\Delta_{N^k}} - \tilde x^k\bigr\|,
\end{equation}
for some $\tau > 0$. Without loss of generality, we can assume that $\bar k$ is sufficiently large to satisfy both the above inequality
and the one of Lemma~\ref{lemma:pg1}.
Therefore, combining Lemma~\ref{lemma:pg1}, \eqref{f_decr} and~\eqref{pg_ineq3}, for $k\ge \hat k$ we can write
\[
\begin{split}
f(\tilde x^k) - f(x^*) & \le \frac L 2 \| \tilde x^k - x^*\|^2 \\[1ex]
                       & \le \frac L 2 \tau^2 \bigl\|P\bigl(\tilde x^k - \nabla f(\tilde x^k)\bigr)_{\Delta_{N^k}} - \tilde x^k\bigr\|^2 \\[1ex]
                       & \le \frac{s L \tau^2}{2 \gamma R \min\{1,s\}^2} \bigl[f(\tilde x^k) - f(x^{k+1})\bigr].
\end{split}
\]
Rearranging the terms and taking into account that $f(\tilde x^k)\le f(x^k)$, we get for all $k \ge \bar k$
\[
 f(x^{k+1}) - f(x^*) \le (1-\rho^{\text{AS-PG}}) \bigl[f(\tilde x^k) - f(x^*)\bigr] \le (1-\rho^{\text{AS-PG}}) \bigl[f(x^k) - f(x^*)\bigr],
\]
where $\displaystyle{\rho^{\text{AS-PG}} = \frac{2 \gamma R \min\{1,s\}^2}{s L \tau^2}}$. 
\end{proof}

\section{Numerical Results}\label{sec:res_simplex}
In this section, we report the numerical experience related to our active-set algorithmic framework.
In the following, we denote by \FW, \AFW\ and \PG\ the standard Frank-Wolfe,
the standard away-step Frank-Wolfe and the standard Projected Gradient method, respectively.
For any feasible point $x$, the search directions for these three algorithms are computed as follows:
\begin{itemize}
\item\textit{standard Frank-Wolfe direction}:
    \[
    d^{\text{FW}} = e_{\hat\imath} - x, \quad \hat\imath \in \argmin_{i=1,\ldots,n} \{\nabla_i f(x)\};
    \]
\item\textit{standard away-step Frank-Wolfe direction}:
    \begin{equation*}
    d^{\text{AFW}} =
    \begin{cases}
    d^{\text{FW}}, \quad & \text{ if } \nabla f(x)^T d^{\text{FW}} \le \nabla f(x)^T d^{\text{A}}, \\
    d^{\text{A}}, \quad & \text{ otherwise},
    \end{cases}
    \end{equation*}
    where
    \[
    d^{\text{A}} = x - e_{\hat\jmath}, \quad  \hat\jmath \in \argmax_{j \colon x_j>0}\{ \nabla_j f(x) \};
    \]
    \smallskip\par
\item\textit{standard projected gradient direction}:
    \[
    d^{\text{PG}} = P\bigl(x -s\nabla f(x)\bigr)_{\Delta} - x,
    \]
    where $s>0$ is a fixed scalar.
\end{itemize}

In the first set of experiments we considered a convex quadratic objective function and the comparisons also included the
the \texttt{P2GP} algorithm~\cite{diserafino:2017}, a gradient-based method specifically devised
for solving quadratic problems with a single linear constraint and bounds on the variables.

In order to have a fair comparison between standard and active-set versions,
the line search used for standard methods is the one described in Algorithm~\ref{alg:ls_as}.

We further denote by \ASFW, \ASAFW\ and \ASPG\, the methods deriving from our algorithmic framework,
where the search direction $d^k$ is computed according to (FW), (AFW) and (PG) rule, respectively.
The main aim of our numerical experience is showing how the use of our active-set strategy can improve the performance of
existing global convergent algorithms.

In our experiments, we set $\delta = 0.5$ and $\gamma = 10^{-4}$ for the Armijo line search,
and $s = 1$ for the computation of $d^k$ when using (PG) rule and the standard projected gradient direction.

In order to calculate the active-set estimate at each iteration of the active-set methods, we need to set the $\eps$ parameter
to a proper value, so that Assumption~\ref{ass:eps} is satisfied.
In general, the value of this parameter cannot be a priori computed.
Following~\cite{cristofari:2017,desantis:2016}, we employ this simple updating rule:
at every iteration $k$, we compute $\tilde x^k$
and, if a sufficient decrease in the objective function is obtained (according to Theorem~\ref{th:estim}),
we accept $\tilde x^k$ and we do not change the value of $\eps$.
Otherwise, we do not accept $\tilde x^k$, we reduce $\eps$ and we estimate the active set again, continuing until
we get a sufficient decrease in the objective function.
The starting value for the $\eps$ parameter is $10^{-1}$ and we set $C=10^{-6}$.

All the codes used in the tests were implemented in Matlab R2014b and the experiments
were ran on an Intel Xeon(R), CPU E5-1650 v2 3.50 GHz.

\subsection{Comparison on Instances from the Chebyshev Center Problem}\label{sec:cheb}
The Chebyshev center problem consists in finding the circle of minimum radius
that encloses all the points in a given finite set $C = \{c_1,\ldots,c_n\}\subset \R^m$
(see, e.g.,~\cite{xu2003solution} and the references therein).
The problem can be formulated as problem~\eqref{simplex_prob},
where $f(x)= x^T {A^T} A x  - \sum_{i=1}^n \|c_i\|^2 x_i$, with
$A = \begin{pmatrix}
      c_1 & & \ldots & & c_n
     \end{pmatrix}
     \in \R^{m\times n}$.
We generated instances with
\begin{itemize}
\item $n \text{ (i.e., cardinality of $C$) } = 2^{15}$;
\item $m \text{ (i.e., samples' dimension) } = 10,\, 100,\, 1000$.
\end{itemize}
For each combination of $n$ and $m$, we considered $10$ problems $p_1,\ldots,p_{10}$ where
we randomly generated the vectors $c_i\in \R^m$, $i=1,\ldots, n$, by
the Matlab's built-in function \textit{rand}.
For each problem $p_h$, $h=1,\ldots,10$, we used the Matlab's command $rng(h)$ to set the seed for the random number generator.
For each problem we fix the starting point to $e_1$ (the same for all methods).

In this set of experiments we first run \ASAFW, stopping it at the first iteration $k$ satisfying
\[
\nabla f(x^k)^T (x-x^k) \ge -10^{-6}, \quad \forall x \in \Delta,
\]
and setting $f_{\text{min}} = f(x^k)$. Then, we used $f_{\text{min}}$ as target value to stop the other algorithms using a certain tolerance
(keep in mind that the problem is convex).
Namely, all the other algorithms,
including \texttt{P2GP}, were stopped at the first iteration $k$ satisfying
\[
f(x^k) \le  f_{\text{min}} + 10^{-6}(1+|f_{\text{min}}|).
\]
We also used a time limit of $3600$ seconds for each method.
In \texttt{P2GP}, all parameters were set to their default values
(except those for termination, that were set in order to stop the algorithm as just described).

In Table~\ref{tab:cheb}, for each instance and each algorithm,
we report the CPU time needed to satisfy the stopping criterion and the
final objective function value found.

\begin{table}
\centering
\caption{Comparison on instances from the Chebyshev center problem. A star indicates that the time limit of $3600$ seconds was reached.}\label{tab:cheb}
\resizebox{15cm}{!} {
\centering
\begin{tabular}{|c c || c | c | c | c | c | c | c || c | c | c | c | c | c | c |}
\hline
\multicolumn{1}{|c}{\multirow{2}*{$m$}} & \multicolumn{1}{c||}{\multirow{2}*{}} & \multicolumn{7}{c||}{CPU time} & \multicolumn{7}{c|}{Obj} \bigstrut[t] \bigstrut[b] \\
& & \FW\ & \ASFW\ & \AFW\ & \ASAFW\ & \PG\ & \ASPG\ & \texttt{P2GP} & \FW\ & \ASFW\ & \AFW\ & \ASAFW\ & \PG\ & \ASPG\ & \texttt{P2GP} \bigstrut[b] \\
\hline
\multirow{10}*{$10$}
& $p_{1}$ & $*$ & $36.22$ & $40.75$ & $0.92$ & $269.21$ & $0.48$ & $146.90$ & $-33.49$ & $-33.49$ & $-33.49$ & $-33.49$ & $-33.49$ & $-33.49$ & $-33.49$ \bigstrut[t] \\
& $p_{2}$ & $*$ & $0.34$ & $12.55$ & $0.37$ & $350.77$ & $0.33$ & $*$ & $-35.78$ & $-35.79$ & $-35.79$ & $-35.79$ & $-35.79$ & $-35.79$ & $-35.77$ \\
& $p_{3}$ & $*$ & $20.47$ & $12.46$ & $0.42$ & $396.47$ & $0.59$ & $652.64$ & $-36.58$ & $-36.58$ & $-36.58$ & $-36.58$ & $-36.58$ & $-36.58$ & $-36.58$ \\
& $p_{4}$ & $*$ & $58.68$ & $14.62$ & $0.39$ & $261.99$ & $0.47$ & $74.20$ & $-34.92$ & $-34.92$ & $-34.92$ & $-34.92$ & $-34.92$ & $-34.92$ & $-34.92$ \\
& $p_{5}$ & $*$ & $0.25$ & $7.62$ & $0.27$ & $389.11$ & $0.34$ & $961.91$ & $-37.56$ & $-37.56$ & $-37.56$ & $-37.56$ & $-37.56$ & $-37.56$ & $-37.56$ \\
& $p_{6}$ & $*$ & $279.68$ & $51.44$ & $1.01$ & $419.54$ & $0.89$ & $*$ & $-35.96$ & $-35.96$ & $-35.96$ & $-35.96$ & $-35.96$ & $-35.96$ & $-35.96$ \\
& $p_{7}$ & $*$ & $0.34$ & $11.48$ & $0.36$ & $216.04$ & $0.41$ & $508.57$ & $-36.01$ & $-36.01$ & $-36.01$ & $-36.01$ & $-36.01$ & $-36.01$ & $-36.01$ \\
& $p_{8}$ & $*$ & $0.54$ & $15.92$ & $0.40$ & $165.88$ & $0.42$ & $3127.75$ & $-35.34$ & $-35.34$ & $-35.34$ & $-35.34$ & $-35.34$ & $-35.34$ & $-35.34$ \\
& $p_{9}$ & $*$ & $16.87$ & $14.63$ & $0.38$ & $336.34$ & $0.35$ & $100.14$ & $-35.60$ & $-35.60$ & $-35.60$ & $-35.60$ & $-35.60$ & $-35.60$ & $-35.60$ \\
& $p_{10}$ & $*$ & $36.90$ & $25.33$ & $0.47$ & $251.71$ & $0.46$ & $*$ & $-33.13$ & $-33.13$ & $-33.13$ & $-33.13$ & $-33.13$ & $-33.13$ & $-33.13$ \bigstrut[b] \\
\hline
\multirow{10}*{$100$}
& $p_{1}$ & $*$ & $332.95$ & $60.12$ & $1.07$ & $1089.13$ & $1.60$ & $22.51$ & $-156.85$ & $-156.86$ & $-156.86$ & $-156.86$ & $-156.86$ & $-156.86$ & $-156.86$ \bigstrut[t] \\
& $p_{2}$ & $*$ & $326.44$ & $51.96$ & $1.64$ & $878.55$ & $2.04$ & $31.94$ & $-153.68$ & $-153.68$ & $-153.68$ & $-153.68$ & $-153.68$ & $-153.68$ & $-153.68$ \\
& $p_{3}$ & $*$ & $74.47$ & $43.54$ & $0.83$ & $137.13$ & $2.45$ & $18.70$ & $-156.05$ & $-156.05$ & $-156.05$ & $-156.05$ & $-156.05$ & $-156.05$ & $-156.05$ \\
& $p_{4}$ & $*$ & $33.49$ & $50.20$ & $1.16$ & $848.55$ & $0.96$ & $22.97$ & $-153.15$ & $-153.16$ & $-153.16$ & $-153.16$ & $-153.16$ & $-153.16$ & $-153.16$ \\
& $p_{5}$ & $*$ & $6.93$ & $43.55$ & $0.97$ & $920.43$ & $2.76$ & $26.49$ & $-152.91$ & $-152.91$ & $-152.91$ & $-152.91$ & $-152.91$ & $-152.91$ & $-152.91$ \\
& $p_{6}$ & $*$ & $326.46$ & $46.81$ & $1.06$ & $908.36$ & $1.39$ & $17.83$ & $-152.75$ & $-152.75$ & $-152.75$ & $-152.75$ & $-152.75$ & $-152.75$ & $-152.75$ \\
& $p_{7}$ & $*$ & $6.68$ & $45.52$ & $1.07$ & $983.68$ & $0.85$ & $*$ & $-151.88$ & $-151.88$ & $-151.88$ & $-151.88$ & $-151.88$ & $-151.88$ & $-151.88$ \\
& $p_{8}$ & $*$ & $317.54$ & $48.29$ & $1.05$ & $1033.65$ & $1.33$ & $217.97$ & $-153.43$ & $-153.43$ & $-153.43$ & $-153.43$ & $-153.43$ & $-153.43$ & $-153.43$ \\
& $p_{9}$ & $*$ & $225.92$ & $51.03$ & $1.06$ & $968.77$ & $1.97$ & $231.58$ & $-152.26$ & $-152.26$ & $-152.26$ & $-152.26$ & $-152.26$ & $-152.26$ & $-152.26$ \\
& $p_{10}$ & $*$ & $52.11$ & $58.23$ & $1.26$ & $967.47$ & $2.25$ & $41.47$ & $-153.34$ & $-153.34$ & $-153.34$ & $-153.34$ & $-153.34$ & $-153.34$ & $-153.34$ \bigstrut[b] \\
\hline
\multirow{10}*{$1000$}
& $p_{1}$ & $*$ & $8.45$ & $151.76$ & $7.38$ & $1770.16$ & $3.76$ & $29.18$ & $-1131.67$ & $-1131.67$ & $-1131.67$ & $-1131.67$ & $-1131.67$ & $-1131.67$ & $-1131.67$ \bigstrut[t] \\
& $p_{2}$ & $*$ & $128.50$ & $140.18$ & $6.62$ & $1374.36$ & $5.30$ & $13.65$ & $-1134.48$ & $-1134.49$ & $-1134.49$ & $-1134.49$ & $-1134.49$ & $-1134.49$ & $-1134.49$ \\
& $p_{3}$ & $*$ & $38.59$ & $165.97$ & $6.30$ & $1114.03$ & $3.43$ & $13.14$ & $-1136.53$ & $-1136.54$ & $-1136.54$ & $-1136.54$ & $-1136.54$ & $-1136.54$ & $-1136.54$ \\
& $p_{4}$ & $*$ & $32.73$ & $155.13$ & $6.55$ & $2009.08$ & $2.98$ & $12.77$ & $-1137.88$ & $-1137.88$ & $-1137.88$ & $-1137.88$ & $-1137.88$ & $-1137.88$ & $-1137.88$ \\
& $p_{5}$ & $*$ & $10.74$ & $153.57$ & $6.37$ & $2284.88$ & $4.88$ & $11.11$ & $-1138.01$ & $-1138.01$ & $-1138.01$ & $-1138.01$ & $-1138.01$ & $-1138.01$ & $-1138.01$ \\
& $p_{6}$ & $*$ & $25.39$ & $179.13$ & $8.07$ & $900.20$ & $3.42$ & $11.34$ & $-1131.69$ & $-1131.70$ & $-1131.70$ & $-1131.70$ & $-1131.70$ & $-1131.70$ & $-1131.70$ \\
& $p_{7}$ & $*$ & $13.82$ & $123.64$ & $6.83$ & $1467.58$ & $8.15$ & $12.76$ & $-1136.24$ & $-1136.26$ & $-1136.26$ & $-1136.26$ & $-1136.26$ & $-1136.26$ & $-1136.26$ \\
& $p_{8}$ & $*$ & $65.76$ & $136.15$ & $5.35$ & $2422.12$ & $3.37$ & $13.49$ & $-1138.78$ & $-1138.79$ & $-1138.79$ & $-1138.79$ & $-1138.79$ & $-1138.79$ & $-1138.79$ \\
& $p_{9}$ & $*$ & $542.44$ & $161.34$ & $6.63$ & $969.59$ & $13.65$ & $11.90$ & $-1136.86$ & $-1136.87$ & $-1136.87$ & $-1136.87$ & $-1136.87$ & $-1136.87$ & $-1136.87$ \\
& $p_{10}$ & $*$ & $37.31$ & $162.30$ & $7.43$ & $1644.57$ & $18.22$ & $118.83$ & $-1135.62$ & $-1135.62$ & $-1135.62$ & $-1135.62$ & $-1135.62$ & $-1135.62$ & $-1135.62$ \bigstrut[b] \\
\hline
\end{tabular}
}
\end{table}

\begin{figure}[h!]
\centering
\includegraphics[scale=0.47, trim = 2.5cm 0.6cm 1.8cm 1cm, clip]{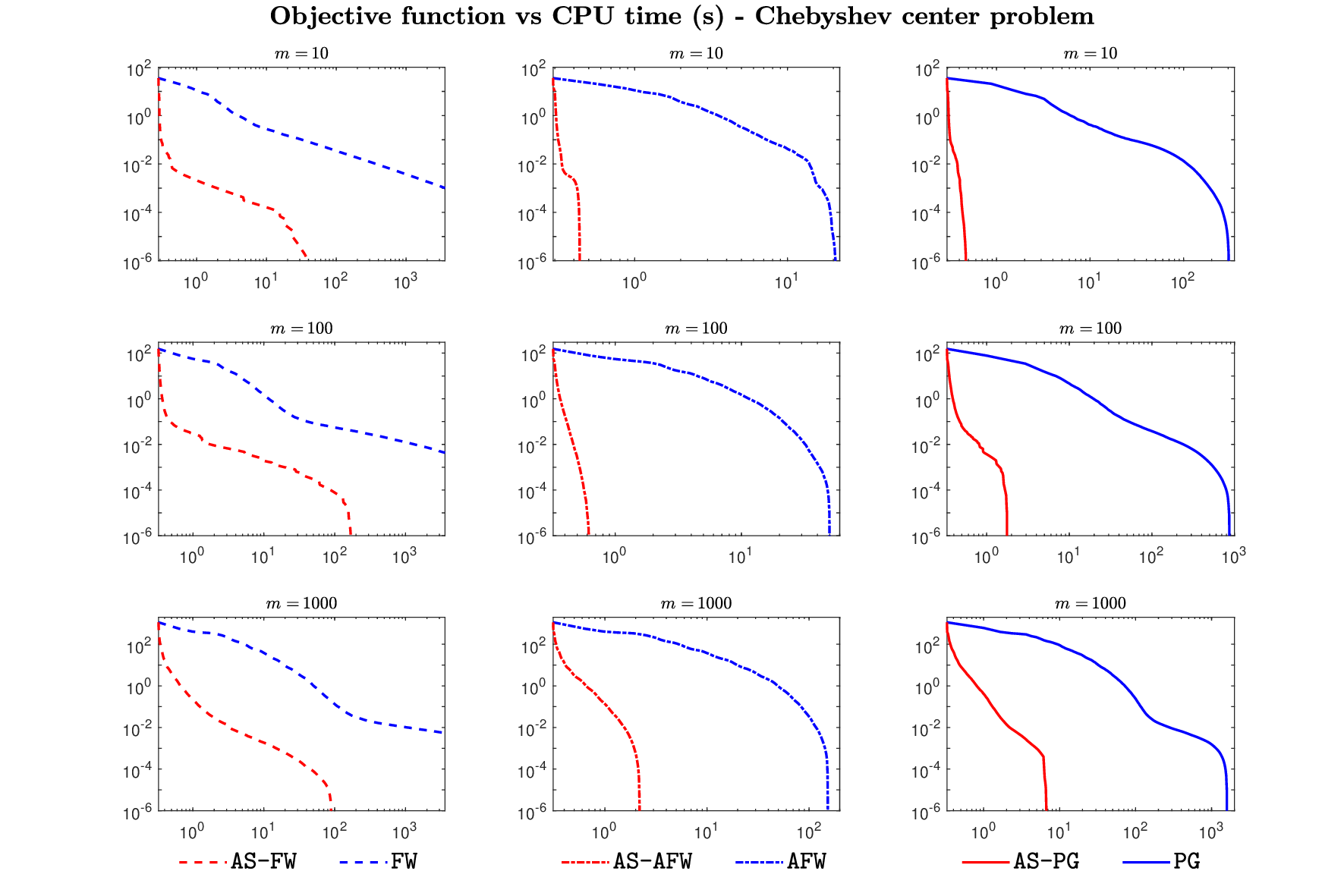}
\caption{Objective function error ($y$ axis) vs CPU time in seconds ($x$ axis).
Comparison between original and active-set algorithms on instances from the Chebyshev center problem.
Both $y$ axis and $x$ axis are in logarithmic scale.}
\label{fig:cheb}
\end{figure}

We first see that the active-set versions are up to two orders of magnitude faster than the standard counterparts.
More in detail, in Figure~\ref{fig:cheb} we report the optimization error for the comparison between \ASFW\ and \FW,
between \ASAFW\ and \AFW, and between \ASPG\ and \PG.
For every $m$, in each plot we report the optimization error $E^k=f(x^k)-f_{\text{min}}$,
averaged over the $10$ runs, versus the computational time.
We can notice that the active-set methods clearly outperform the standard algorithms.

\begin{figure}[h!]
\centering
\includegraphics[scale=0.47, trim = 2.5cm 12.85cm 1.8cm 1cm, clip]{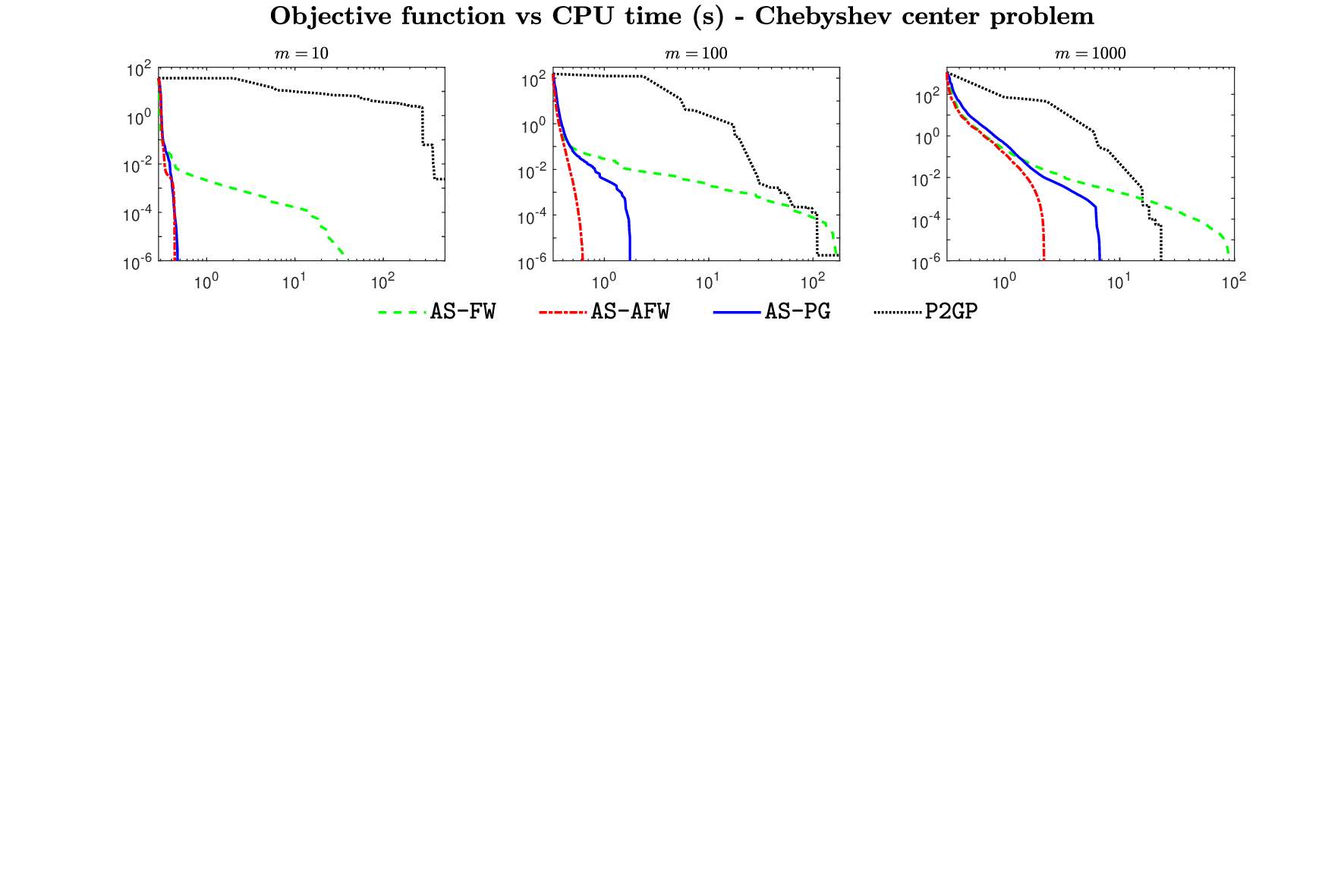}
\caption{Objective function error ($y$ axis) vs CPU time in seconds ($x$ axis).
Comparison among active-set algorithms and \texttt{P2GP} on instances from the Chebyshev center problem.
Both $y$ axis and $x$ axis are in logarithmic scale.}
\label{fig:chebyshev_all_solvers}
\end{figure}

For what concerns the comparisons among our active-set methods and \texttt{P2GP}, from Table~\ref{tab:cheb} we see that
both \ASAFW\ and \ASPG\ are much faster than \texttt{P2GP} for all the considered values of $m$,
while \ASFW\ is faster than \texttt{P2GP} for $m=10$, slower for $m=1000$ and is comparable with \texttt{P2GP} for $m=100$
(in the sense that on average they take an amount of time of the same order of magnitude).
For these comparisons, the optimization error $E^k=f(x^k)-f_{\text{min}}$, averaged over the $10$ runs, is reported in Figure~\ref{fig:chebyshev_all_solvers}.

\subsection{Comparison on Instances from the Symmetric Eigenvalue Complementarity Problem}\label{sec:eicp}
In our experiments, we further consider a class of non-convex non-quadratic problems, namely instances from the
symmetric Eigenvalue Complementarity Problem (EiCP).
Given the symmetric matrices $A \in \R^{n\times n}$ and $B\in \R^{n\times n}$, $B\succ 0$, the symmetric EiCP consists
of finding a scalar $\lambda$ and a vector $x \in \R^n\setminus\{0\}$ such that
\begin{align*}
& w = (\lambda B - A)x, \\
& w \ge 0, \; x\ge 0, \\
& w^T x = 0.
\end{align*}
It can be proven that the symmetric EiCP reduces to computing a stationary point of problem~\eqref{simplex_prob}
where $f(x)= -\frac{x^T A x}{x^T B x}$~\cite{iusem2019splitting}.

Taking into account the numerical experience reported in~\cite{iusem2019splitting},
we considered instances with $B$ equal to the identity matrix and $A\in \R^{n\times n}$ a negative definite matrix.
More in detail, we considered $10$ problems $p_1,\ldots,p_{10}$ where $-A$ was generated as proposed in \cite{more1989algorithms}, with $n = 2^{15}$.
Starting from a vector $y \in \R^n$ with components randomly uniformly distributed in $[-1,1]$ (obtained by the Matlab's built-in function \textit{rand}),
we set the matrix $Y$ as
\[
Y = I -2\frac{yy^T}{\|y\|^2},
\]
where $I$ is the identity matrix and, for each problem $p_h$, $h=1,\ldots,10$, we used the Matlab's command $rng(h)$ to set the seed for the random number generator.
The matrix $A$ was then generated as $-YDY$, where $D$ is a diagonal matrix
with entries $D_{ii} = e^{\frac{i-1}{n-1}}$, $i=1,\ldots,n$.
For each problem, we randomly generated a feasible starting point (the same for all methods).

All algorithms were stopped at the first iteration $k$ satisfying
\[
\nabla f(x^k)^T (x-x^k) \ge -10^{-4}, \quad \forall x \in \Delta.
\]
We also considered a time limit of $3600$ seconds (as for the experiments of the previous subsection),
but it was never reached by any considered method.

\begin{figure}[h!]
\centering
\includegraphics[scale=0.47, trim = 2.5cm 6.6cm 1.8cm 1cm, clip]{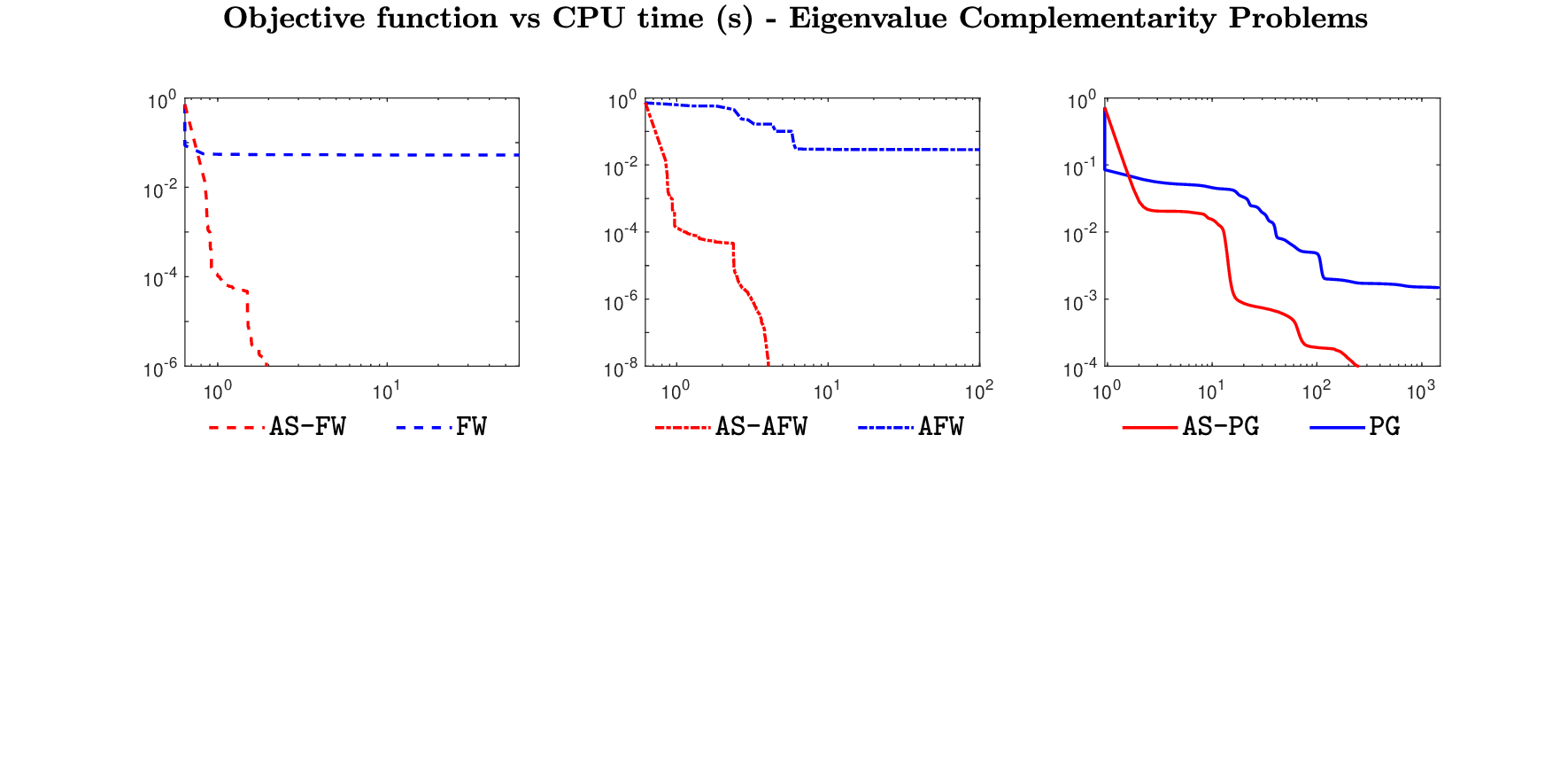}
\caption{Objective function error ($y$ axis) vs CPU time in seconds ($x$ axis).
Comparison between original and active-set algorithms on instances from the symmetric eigenvalue complementarity problem.
Both $y$ axis and $x$ axis are in logarithmic scale.}
\label{fig:eicp}
\end{figure}

In Table~\ref{tab:eicp}, for each instance and each algorithm,
we report the CPU time needed to satisfy the stopping criterion and the
final objective function value found.
In Figure~\ref{fig:eicp}, we report the optimization error for the comparison between \ASFW\ and \FW,
between \ASAFW\ and \AFW, and between \ASPG\ and \PG. In particular, in each plot we report the optimization error $E^k=f(x^k)-f_{\text{min}}$,
averaged over the $10$ runs, versus the computational time, where $f_{\text{min}}$ is the smallest objective value
found in any comparison between the two considered algorithms.

\begin{table}
\centering
\caption{Comparison on instances from the symmetric Eigenvalue Complementarity Problem.}\label{tab:eicp}
\resizebox{15cm}{!} {
\centering
\begin{tabular}{| c || c | c | c | c | c | c || c | c | c | c | c | c |}
\hline
\multicolumn{1}{|c||}{\multirow{2}*{}} & \multicolumn{6}{c||}{CPU time} & \multicolumn{6}{c|}{Obj} \bigstrut[t] \bigstrut[b] \\
& \FW\ & \ASFW\ & \AFW\ & \ASAFW\ & \PG\ & \ASPG\ & \FW\ & \ASFW\ & \AFW\ & \ASAFW\ & \PG\ & \ASPG\ \bigstrut[b] \\
\hline
$p_{1}$ & $56.57$ & $1.62$ & $35.73$ & $1.15$ & $3148.11$ & $40.47$ & $1.00$ & $1.00$ & $1.00$ & $1.00$ & $1.00$ & $1.00$ \bigstrut[t] \\
$p_{2}$ & $31.41$ & $1.15$ & $140.99$ & $2.77$ & $258.68$ & $177.18$ & $1.00$ & $1.00$ & $1.00$ & $1.00$ & $1.00$ & $1.00$ \\
$p_{3}$ & $265.40$ & $11.85$ & $345.04$ & $15.75$ & $1602.88$ & $1557.03$ & $1.00$ & $1.00$ & $1.00$ & $1.00$ & $1.00$ & $1.00$ \\
$p_{4}$ & $0.66$ & $1.01$ & $52.39$ & $1.21$ & $2056.99$ & $149.71$ & $1.24$ & $1.00$ & $1.00$ & $1.00$ & $1.00$ & $1.00$ \\
$p_{5}$ & $0.69$ & $0.82$ & $5.13$ & $0.81$ & $240.35$ & $254.33$ & $1.07$ & $1.07$ & $1.07$ & $1.07$ & $1.00$ & $1.00$ \\
$p_{6}$ & $0.63$ & $1.13$ & $0.97$ & $1.72$ & $1531.66$ & $654.37$ & $1.11$ & $1.00$ & $1.11$ & $1.00$ & $1.00$ & $1.00$ \\
$p_{7}$ & $3.08$ & $0.83$ & $6.91$ & $0.80$ & $1224.86$ & $622.67$ & $1.00$ & $1.00$ & $1.00$ & $1.00$ & $1.00$ & $1.00$ \\
$p_{8}$ & $51.59$ & $1.35$ & $121.98$ & $2.17$ & $835.16$ & $35.57$ & $1.00$ & $1.00$ & $1.00$ & $1.00$ & $1.00$ & $1.00$ \\
$p_{9}$ & $217.15$ & $4.86$ & $376.64$ & $7.86$ & $3444.62$ & $41.61$ & $1.00$ & $1.00$ & $1.00$ & $1.00$ & $1.00$ & $1.00$ \\
$p_{10}$ & $0.66$ & $2.72$ & $3.67$ & $6.51$ & $5.07$ & $135.14$ & $1.18$ & $1.00$ & $1.18$ & $1.00$ & $1.01$ & $1.00$ \bigstrut[b] \\
\hline
\end{tabular}
}
\end{table}

Also for these problems, we can notice that the active-set framework clearly
outperforms on average the original algorithms used in the comparison.
In particular, the active-set versions are on average able to find a better solution in a smaller amount of time.

\section{Conclusions}\label{sec:conclusions_simplex}
In this paper, we focused on minimization problems over the simplex
and described an active-set algorithmic framework.
The active-set strategy we adopted here does not only focus on the zero variables and keep them fixed, but rather tries
to quickly identify as many active variables as possible (including nonzero variables)
at a given point. Furthermore, it suitably reduces the objective function
(when setting to zero those variables estimated active), while guaranteeing
feasibility. This last feature,  together  with the use of
active-set gradient related directions and an Armijo line search, allowed us to prove global
convergence of the framework.
We further described three different types of active-set gradient related directions  and
proved linear convergence rate when using those directions in the algorithm.
Our numerical experience on sparse optimization problems highlighted
the efficiency of our new method when dealing with both non-convex and convex instances.

\bibliography{as_simplex}

\end{document}